\documentclass[preprint,12pt]{elsarticle}

\usepackage[pdftex]{hyperref}
\biboptions{sort&compress,numbers,square}

\usepackage{amsfonts}
\usepackage{graphicx}
\usepackage{epsfig}
\usepackage{geometry}
\usepackage{tikz}
\usetikzlibrary{trees,positioning,shapes}
\usepackage{subcaption}
\usepackage{amsmath}
\usepackage{amsthm}

\usepackage{xcolor, colortbl}

\usepackage{float}
\usepackage[hypcap]{caption}
\usepackage{subcaption}
\usepackage{algorithm2e}
\usepackage{array}

\newcommand{\inblue}[1]{\textcolor{blue}{#1}}

\newtheorem{definition}{Definition}

\newtheorem{remark}{Remark}

\usepackage{bm}
\usepackage{mathtools}

\usepackage{rotating}

\begin{document}

\begin{frontmatter}

\title{Quasi-optimal $hp$-finite element refinements towards singularities via deep neural network prediction}

\author{Tomasz S\l{}u\.zalec$^{(1)}$, Rafa\l{} Grzeszczuk$^{(2)}$,  Sergio Rojas$^{(3)}$, Witold Dzwinel $^{(2)}$,
\\ Maciej Paszy\'{n}ski$^{(2)}$}

\address{$^{(1)}$Jagiellonian University, Krak\'ow, Poland \\ $^{(2)}$Institute of Computer Science, \\ AGH University of Science and Technology,
Krak\'{o}w, Poland \\
e-mail: maciej.paszynski@agh.edu.pl \\ $^{(3)}$ Instituto de matem\'aticas, \\
Pontificia Universidad Cat\'olica de Valpara\'iso, Valpara\'iso, Chile\\
e-mail: sergio.rojas.h@pucv.cl}

\begin{abstract}
We show how to construct the deep neural network (DNN) expert to predict quasi-optimal $hp$-refinements for a given computational problem. The main idea is to train the DNN expert during executing the self-adaptive $hp$-finite element method ($hp$-FEM) algorithm and use it later to predict further $hp$ refinements.
For the training, we use a two-grid paradigm self-adaptive $hp$-FEM algorithm. It employs the fine mesh to provide the optimal $hp$ refinements for coarse mesh elements.
We aim to construct the DNN expert to identify quasi-optimal $hp$ refinements of the coarse mesh elements.
During the training phase, we use the direct solver to obtain the solution for the fine mesh to guide the optimal refinements over the coarse mesh element. 
After training, we turn off the self-adaptive $hp$-FEM algorithm and continue with quasi-optimal refinements as proposed by the DNN expert trained. We test our method on three-dimensional Fichera and two-dimensional L-shaped domain problems. We verify the convergence of the numerical accuracy with respect to the mesh size. We show that the exponential convergence delivered by the self-adaptive $hp$-FEM can be preserved if we continue refinements with a properly trained DNN expert. Thus, in this paper, we show that from the self-adaptive $hp$-FEM it is possible to train the DNN expert the location of the singularities, and continue with the selection of the quasi-optimal $hp$ refinements, preserving the exponential convergence of the method.
\end{abstract}
	

\end{frontmatter}

\section{Introduction}

The application of the neural networks in adaptive finite element method usually supports the generation of the adaptive mesh \cite{WDz2} or mesh partitioning algorithms \cite{WDz1}.
In this paper, we focus on $hp$ adaptive finite element method. It can approximate a given problem with prescribed accuracy using a minimal number of basis functions.
It employs a mesh to approximate a given problem using a linear combination of basis functions spanning the computational mesh.
It improves approximation quality by generating new basis functions for a better approximation.
The new basis functions are obtained by breaking selected elements into smaller elements and spanning smaller basis functions there (this process is called $h$ refinement). We may also obtain them by adding new higher-order basis functions over existing elements (this process is called the $p$ refinement). The $hp$ refinement is the combination of both.
It has been proven theoretically and experimentally \cite{Babuska1,Babuska2,Mitchell}, that the $hp$ adaptive algorithm can deliver the best possible exponential convergence of the numerical error with respect to the mesh size.

There are several implementations of the data structures supporting $hp$-refinements \cite{Book1,Book2,Rachowicz2,Bank,Rank,Multilevel}.
In this paper, we focus on the data structure proposed by Demkowicz \cite{Book2,Rachowicz2}, supporting hanging nodes and the 1-irregularity rule.
This code also contains the implementation of the fully automatic self-adaptive $hp$ finite element method algorithm (self-adaptive $hp$-FEM), described in \cite{Rachowicz,Rachowicz2}.
It employs the two-grid paradigm.
It compares numerical solutions on the two grids, the actual (the coarse) and the reference (the fine) grid. Each finite element from the coarse mesh, it considers different refinement strategies. For each strategy, it computes the error decrease rate.
It estimates how much error decrease we gain by selecting this local element refinement, divided by how many degrees of freedom we have to invest to get there. For evaluating the error decrease rate, it employs the projections from the fine mesh local element solution into the considered refinement strategy of the coarse mesh element. This strategy selects the optimal refinements for particular elements. It delivers the exponential convergence, which has been verified experimentally for a large class of problems \cite{Book1,Book2}.

However, this implementation of the self-adaptive $hp$-FEM algorithm is expensive. It requires the solutions to the computational problem over the coarse and fine mesh in several iterations. In particular, the fine mesh problem is usually one order of magnitude larger than the actual coarse mesh problem in two dimensions and two orders of magnitude larger in 3D.
The computational grids are highly non-uniform; they contain elongated elements and different polynomial orders of approximations. Thus, the iterative solvers deliver convergence problems \cite{TwoGridSolver,Saad,static}, and the computationally expensive \cite{Singularities1,Singularities2} direct solvers \cite{DS1,DS2} augmented with proper ordering algorithms \cite{c23,c24,c25} are employed to solve the coarse and the fine mesh problems. Additionally, selecting the optimal refinements involves computations of local projections and multiple possibilities to consider over the coarse mesh elements.

Here we propose an artificial expert, Deep Neural Network (DNN), selecting refinements for particular elements of the computational mesh. The DNN expert focus on the key problem in the implementation of the self-adaptive $hp$-FEM, namely the selection of the optimal refinements for the mesh elements. We train the DNN expert during the execution of the self-adaptive $hp$-FEM algorithm.
We collect the dataset for all the decisions made by the self-adaptive $hp$-FEM. The input to the DNN is element data, such as its location, sizes, and actual polynomial orders of approximations. The output from the DNN is the decision about optimal refinements. We train the DNN to make quasi-optimal decisions about optimal refinements.
When the DNN is ready to make quasi-optimal decisions, we turn off the computationally intense self-adaptive $hp$-FEM, and we ask the DNN-driven $hp$-FEM to predict the optimal refinements. Our DNN-based expert can learn the locations of the singularities where the $h$ refinements are needed, as well as the patterns for performing the $p$ refinements.

The DNN has already been applied for $hp$-adaptive finite element method in two-dimensions \cite{DNNhp}. The experiments presented there concern the case of the DNN training to learn the optimal $hp$-refinements for a given class of boundary-value problems to apply the trained DNN for new problems from the same class.
Our paper differs from \cite{DNNhp} in three aspects. First, we generalize the DNN-driven $hp$-FEM algorithm into three dimensions. Second, we train the DNN expert to propose further quasi-optimal $hp$ refinements for a given computational problem when the self-adaptive algorithm runs out of resources or becomes computationally too intense. Third, in \cite{DNNhp}, the DNN is trained to make decisions about local element refinement by using the input data, including the local solution. We can only use the DNN from \cite{DNNhp} after we call a solver. This makes it impossible to use that DNN to guide further refinements without calling the solver.
In this paper, we train the DNN for the two- and 3D problem, using as input the element coordinates and the actual polynomial orders of approximation over the element. In this way, our DNN-based expert learns the locations of singularities, the geometric patterns for $h$ refinements, and the patterns for modification of the polynomial orders of approximation.
Thus, we propose a DNN expert that allows continuing with $hp$ refinements beyond the reach of the solver.

The DNN-based expert training methodology proposed in this paper can guide refinements in different computational codes, including the $hp$ codes developed by the group of Leszek Demkowicz \cite{Book1,Book2}, or the codes using a hierarchy of nested meshed for performing the refinements \cite{Painless,Rank,Multilevel,Rank1,Rank2}. The generalization of the idea presented in this paper to different basis functions and multiple hanging nodes requires modification of the input and/or output of the DNN. The DNN, once trained, proposes quasi-optimal $hp$ refinements for finite elements from the actual computational mesh.

The structure of the paper is the following. Section 2 is the Preliminaries where we introduce our model problem, some basic definitions and the self-adaptive $hp$-FEM algorithm. Section 3 introduces the architecture of the neural network and the DNN-driven $hp$-FEM algorithm.  Section 4 presents the numerical results for the model L-shape domain problem and the model Fichera problem defined as the superposition of several model L-shaped domain problems.  We conclude the paper in Section 5.

\section{Preliminaries}
\subsection{Model problem and variational formulations}\label{ch2:model}
Let $\Omega \subset \mathbb{R}^d$,  with $d=2,3$,  be a bounded, open and simply connected polyhedron, and denote by $\partial \Omega$ its boundary assumed to be split into two disjoint open and not empty subsets $\Gamma_D,  \Gamma_N$.  That is, $ \Gamma_D \cup \Gamma_N = \emptyset$ and $\overline{\partial \Omega} = \overline{\Gamma_D \cup \Gamma_N}$.  Using the standard notations for Hilbert spaces,  for a given source term $g \in L^2\left(\Gamma_N \right)$,  we look for an approximation of the solution $u$ satisfying:
\begin{equation}\label{eq:laplace}
-\Delta u = 0 \text{ in } \Omega, \qquad 
u = 0 \text{ on } \Gamma_D, \qquad 
\frac{\partial u}{\partial n} =g \text{ on } \Gamma_N.
\end{equation}
Introducing the Hilbert space $V :=\left\{ v \in H^1(\Omega) \, : \,  v = 0 \text{ on } \Gamma_D \right\}$,  problem \eqref{eq:laplace} admits the following well-posed continuous variational formulation \cite{Book2}:
\begin{equation}\label{eq:cvf}
\displaystyle \text{Find } u \in V \, : \, b(u,v) := \int_\Omega \nabla u \cdot \nabla v \, dx = \int_{\Gamma_N} g v \, dS =: l(v),  \, \forall \, v \in V.
\end{equation}
%
%
The discrete formulation
\begin{equation}\label{eq:cvf}
\displaystyle \text{Find } u \in V_h \, : \, b(u_h,v_h) := \int_\Omega \nabla u_h \cdot \nabla v_h \, dx = \int_{\Gamma_N} g v_h \, dS =: l(v_h),  \, \forall \, v_h \in V_h.
\end{equation}
The discrete space $V_h$ is constructed with hierarchical basis functions obtained by glueing together the element shape functions.
The shape functions are constructed as tensor products of 1D hierarchical-shape functions


\subsection{Self-adaptive $hp$ finite element method}

In this Section, we briefly discuss the self-adaptive $hp$-FEM algorithm,  following the ideas described by Demkowicz et al. \cite{Book1,Book2}.  We,  however, propose a slightly different algorithm,  with the differences explained inside this Section.  \\
%

The self-adaptive $hp$-FEM algorithm is based on a two-grid paradigm.  A computation over a coarse mesh,  and a computation over a fine mesh.

\begin{definition} 
    The \textit{initial coarse mesh} is obtained by partitioning the domain $\Omega$ into a finite set $(K, X(K), \Pi_p) \in T_{hp}$ of \textit{hp} finite elements and selecting arbitrary polynomial orders of approximation. \\
\end{definition}
\begin{definition} 
    \textit{Coarse mesh problem}: Find $\{u_{hp}^i\}_{i=1}^{N_{hp}}$ coefficients (\textit{dofs}) of approximate solution $V \supset V_{hp} \ni u_{hp} = \sum_{i=1}^{N_{hp}}u_{hp}^ie^i_{hp}$ fulfilling (\ref{eq:Fichera}).
\end{definition}
\begin{definition} 
    The \textit{coarse mesh approximation space} is defined as
$$V_{hp} = \text{span}\{e_{hp}^j: \forall K \in T_{hp}|_K, \forall \phi_k \in X(K), \exists! e_{hp}^i:e_{hp}^i|_K = \phi_k \} $$
    where $e_{hp}^i$ is a global basis function (element basis of $V_{hp}$), 
$\phi_k$ is a shape function and $(k, K) \rightarrow i(k_1,K)$ is the mapping 
over the coarse mesh assigning global number $i(k,K)$ of \textit{dofs} (basis functions) related with shape function $k$ from element $K$ 
\end{definition} 

\begin{remark} The approximation space $V_{hp} \subset V$ with basis $\{e_{hp}^i\}_{i=1}^{N_{hp}}$ is constructed by gluing together element-local-shape functions.
\end{remark}

%
%

\begin{definition} 
    The \textit{fine mesh} is obtained by breaking each element from the coarse mesh  $(K, X(K), \Pi_p) \in T_{\frac{h}{2}, p+1}$ into 8 elements (in 3D) and increasing the polynomial orders of approximation by one.
\end{definition}

\begin{definition}     \textit{Fine mesh problem}: Find $\{u_{\frac{h}{2},p+1}^i\}_{i=1}^{N_{{\frac{h}{2},p+1}}}$ coefficients (\textit{dofs}) of approximate solution $V \supset V_{{\frac{h}{2},p+1}} \ni u_{{\frac{h}{2},p+1}} = \sum_{i=1}^{N_{{\frac{h}{2},p+1}}}u_{{\frac{h}{2},p+1}}^ie^i_{{\frac{h}{2},p+1}}$ fulfilling (\ref{eq:Fichera}).
    
\end{definition}

\begin{definition} 
The \textit{fine mesh approximation space} is defined as
     \begin{align}
    V_{\frac{h}{2},p+1} =\text{span}
    \{e_{\frac{h}{2},p+1}^j: \forall K \in T_{\frac{h}{2},p+1}|_K, \forall \phi_k \in X(K), \exists! e_{\frac{h}{2},p+1}^i:e_{\frac{h}{2},p+1}^i|_K = \phi_k \} \nonumber
    \end{align}
    where $e_{\frac{h}{2},p+1}^i$ is a basis function (element basis of $V_{\frac{h}{2},p+1}$), $\phi_k$ is a shape function, $(k,K)\rightarrow i(k_1,K)$ is the mapping over the fine mesh assigning global number $i(k,K)$ of \textit{dofs} (basis function) related to shape function $k$ from element $K$.
\end{definition}

    \begin{figure}
        \centering
        \includegraphics[width=0.24\textwidth]{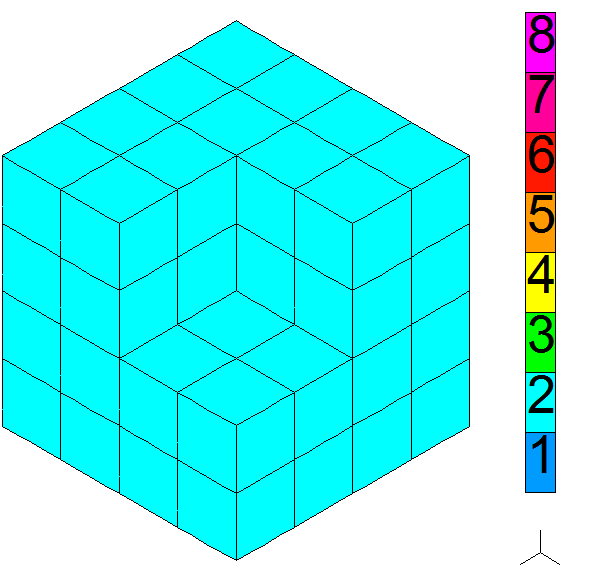}\includegraphics[width=0.24\textwidth]{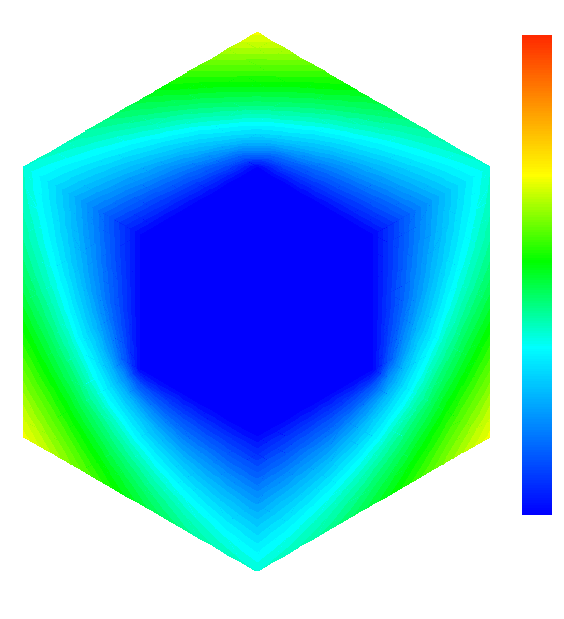} : \includegraphics[width=0.24\textwidth]{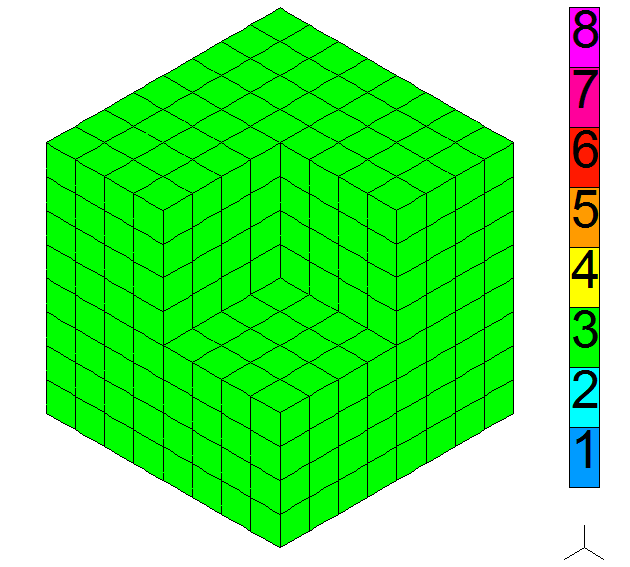}\includegraphics[width=0.24\textwidth]{solution.png} \\       
\noindent\rule{8cm}{0.4pt} \\
        \includegraphics[width=0.25\textwidth]{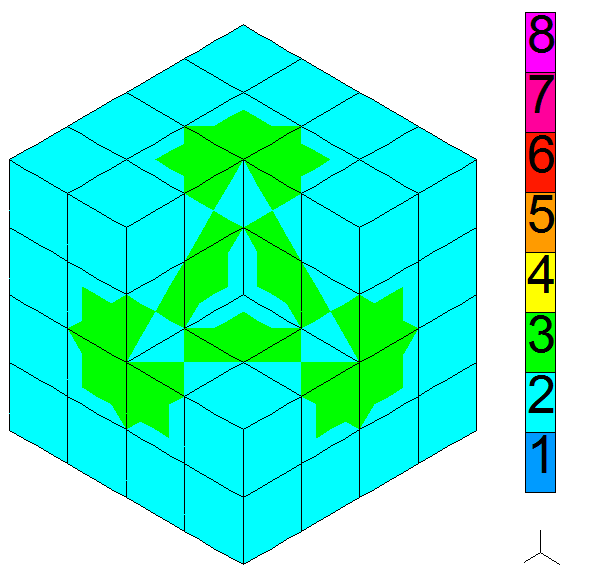}        
\caption{Self-adaptive $hp$ finite element method algorithm}
        
    \end{figure}

    \begin{algorithm}[H]
       
        \caption{Self-adaptive \textit{hp}-FEM algorithm}
        \KwIn{Initial mesh $T_{init}$, PDE, boundary conditions, accuracy}
        \KwOut{Optimal mesh}
        1 $T_{hp}=T_{init}$ (coarse mesh = initial mesh) \\
        2 Solve the coarse mesh problem using direct solver to obtain $u_{hp}$ \\
        3 Generate fine mesh $T_{\frac{h}{2},p+1}$ by breaking each element of the coarse mesh into 8 elements and increasing polynomial orders by 1 \\
        4 Solve the fine mesh problem using direct solver to obtain $u_{\frac{h}{2},p+1}$\\       
        5 Compute $max\_rel\_error = max_{K \in T_{hp}} 100 \times \frac{\| u_\frac{h}{2},p+1 - u_{hp}\|}{\| u_\frac{h}{2},p+1\|}$ \\
        6 \If{maximum relative error $max\_rel\_error <$  accuracy}{
            \Return{$T_{\frac{h}{2}p+1}$ (coarse mesh)}    
          }
        7 Select optimal refinements $\{V_{opt}^K\}_{K\in T_{hp}}$ for every element $K$ from the coarse mesh (\textbf{$\{V^K_{opt}\}_{K \in T_h} = $  Algorithm 2 $(T_{hp},u_{hp},T_{\frac{h}{2},p+1},u_{\frac{h}{2},p+1})$})\\
        8 Perform all \textit{hp} refinements from $\{V_{opt}^K\}_{K\in T_{hp}}$ to obtain $T_{opt}$\\
        9 $T_{hp}=T_{opt}$ (coarse mesh = optimal mesh)\\
        10 \textbf{goto 2}
    \end{algorithm}

Having the self-adaptive $hp$-FEM algorithm defined in Algorithm 1, we focus on Algorithm 2, selecting the optimal $hp$ refinement for a given coarse mesh element.

\begin{figure}
        \centering
        \includegraphics[width=0.6\textwidth]{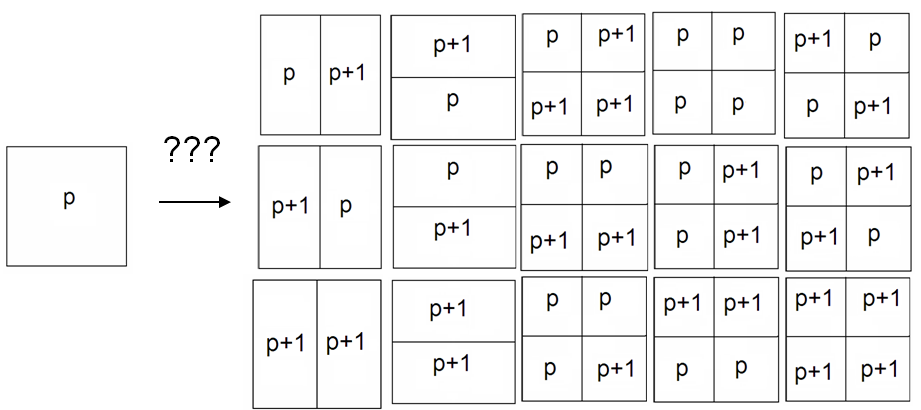}        
       \caption{Selection of optimal refinements in DNN-driven $hp$-FEM code \\ (this simplified picture does not include different orders on edges)}
\label{fig:selection}
    \end{figure}

\begin{figure}
        \centering
        \includegraphics[width=1.0\textwidth]{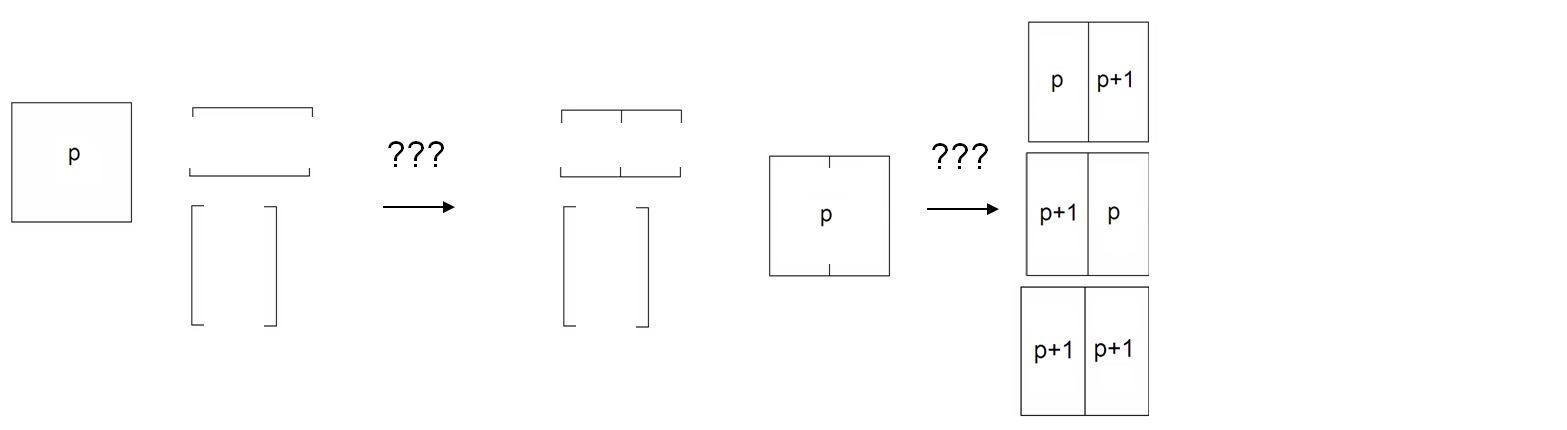}        
       \caption{Selection of optimal refinements in $hp$-FEM code \\ (this simplified picture does not include different orders on edges)}
\label{fig:restriction}
    \end{figure}

Our Algorithm 2, however, differs from the one described by Demkowicz \cite{Book2} in the following aspect. We select optimal refinements for interiors of finite elements (see Figure \ref{fig:selection} for the two-dimensional case), and we adjust the faces and edges accordingly. In the original algorithm described in \cite{Book2}, the optimal refinements are first selected for finite element edges, then for faces, and finally for the interiors. 
The recursive structure of the original algorithm as proposed by \cite{Book2} enables reduction of the computational cost since the selection of the optimal refinements for edges restricts the number of possibilities for refinements of faces, see Figure \ref{fig:restriction}. Similarly, the selection of optimal refinements for faces restricts the number of possible refinements for the interiors.

Our motivation to replace this original three-step algorithm proposed by Demkowicz \cite{Book2} by the one-step Algorithm 2 is the following. We plan to replace the Algorithm 2 by DNN. Once trained, the DNN can provide the quasi-optimal refinement for the interior of the element in a linear computational cost (very fast).

In the one-step Algorithm 2, the selection of the optimal refinements is performed by computing the \emph{error decrease rate} $rate(w)$ for all considered refinement strategies. The $rate(w)$ is defined as the gain in the reduction of the approximation error divided by the loss defined as the number of degrees of freedom added.
The reduction of the approximation error is estimated as $\Big|u_{\frac{h}{2},p+1}-u_{hp}\Big|_{H^1(K)}-\Big|u_{\frac{h}{2},p+1}-w\Big|_{H^1(K)}$, where $\Big|u_{\frac{h}{2},p+1}-u_{hp}\Big|_{H^1(K)}$ denotes the relative error of the coarse mesh element solution $u_{hp}$ with respect to the fine mesh solution $u_{h/2,p+1}$ computed in $H^1$ norm over the element $K$, and $\Big|u_{\frac{h}{2},p+1}-w\Big|_{H^1(K)}$ denotes the relative error of the solution $w$ corresponding to the proposed refinement strategy with respect to the the fine mesh solution $u_{h/2,p+1}$. The difference is the error decrease (how much the error go down if we switch from $u_{hp}$ into $w$ be performing the proposed refinement of the coarse mesh element solution.
For each considered refinement strategy, the solution $w$ corresponding to the refined coarse mesh element is obtained by the projection from the fine mesh solution $u_{h/2,p+1}$.
The procedure is called projection-based interpolation, and it is described in \cite{PBI}.

\begin{algorithm}[H]       
        \caption{Selection of optimal refinements over \textit{K}}
        \SetAlgoLined
        \KwIn{Coarse mesh elements $T_{hp}$, coarse mesh solution $u_{hp} \in V_{hp}$, fine mesh elements $T_{\frac{h}{2},p+1}$, fine mesh solution $u_{\frac{h}{2},p+1} \in V_{\frac{h}{2},p+1}$}
        \KwOut{Optimal refinement $\{ V_{opt}^K \}_{K\in T_{hp}}$ computed for each element $K$}
        \For{coarse mesh elements $K \in T_{hp}$}{
            $rate_{max} = 0$\\
            \For{approximation space $V_{opt} \in K$}{
            Compute the projection based interpolant $w|_K$ of $u_{\frac{h}{2},p+1}|_K$\\
            Compute the error decrease rate 
            $rate(w) = \frac{\Big|u_{\frac{h}{2},p+1}-u_{hp}\Big|_{H^1(K)}-\Big|u_{\frac{h}{2},p+1}-w\Big|_{H^1(K)}}{\Delta\text{nrdof}(V_{hp}, V_{opt}^K, K)}$ \\
            \If{$rate(w) > rate_{max}$}{
                $rate_{max} = rate(w)$\\
                Select $V_{opt}^K$ corresponding to $rate_{max}$ as the optimal refinement for element $K$            
            }
            }
        }
    Select orders on faces = MIN (orders from neighboring interiors) \\
    Select orders on edges = MIN (orders from neighboring faces)
    \end{algorithm}

The fine mesh elements $T_{\frac{h}{2},p+1}$ are employed to localize the 8 elements corresponding to given coarse mesh element $K$, and to compute the projections of the fine mesh solution $u_{\frac{h}{2},p+1}$ into the interpolant $w|_K$.

The following definition formally explains that the optimal refinement for a given element $K$ is the refinement that provides the maximum error decrease rate.

\begin{definition} 
    Let $V_{hp} \subset V_{\frac{h}{2}, p+1} \subset V$ be the coarse and fine mesh approximation spaces. Let $T_{hp}$ represents the coarse mesh elements. \\ 
$u_{hp} \in V_{hp}$ and $u_{\frac{h}{2},p+1} \in V_{\frac{h}{2},p+1}$ are coarse / fine mesh problem solutions. \\
The approximation space $V_{opt}^K$ is called the optimal approximation space over an element $K \in T_{hp}$, if the projection based interpolant $w_{opt}$ of $u_{\frac{h}{2},p+1} \in V_{\frac{h}{2},p+1}$ into $V_{opt}^K$ over element $K$ realizes the following maximum
    \begin{align} \nonumber
        \frac{\Big|u_{\frac{h}{2},p+1}-u_{hp}\Big|_{H^1(K)}-\Big|u_{\frac{h}{2},p+1}-w_{opt}\Big|_{H^1(K)}}{\Delta\text{nrdof}(V_{hp}, V_{opt}^K, K)} &= \\ \nonumber
        =\max_{V_{hp} \subseteq V_{w} \subseteq V_{\frac{h}{2}, p+1}} \frac{\Big|u_{\frac{h}{2},p+1}-u_{hp}\Big|_{H^1(K)}-\Big|u_{\frac{h}{2},p+1}-w\Big|_{H^1(K)}}{\Delta\text{nrdof}(V_{hp}, V_{opt}^K, K)}
    \end{align}
    where $w$ is the projection-based interpolant of $u_{\frac{h}{2},p+1} \in V_{\frac{h}{2},p+1}$ into $V_w$ over element $K$, and $\Delta nrdof(V,X,K) = dim V\Big|_K - dim X\Big|_K$
\end{definition}

Note that the selection of the optimal refinement in the Algorithm 2 is based on the coarse and fine mesh solutions. We plan to train the DNN to select quasi-optimal refinements based on the location of elements and actual polynomial orders of approximation, to continue with the refinement once the solver runs out of resources.

\section{Deep neural network-driven $hp$-adaptive finite element method algorithm}

We propose the following architecture of the DNN that learns the quasi-optimal refinements.

    \begin{figure}
        \centering
        \includegraphics[width=0.9\textwidth]{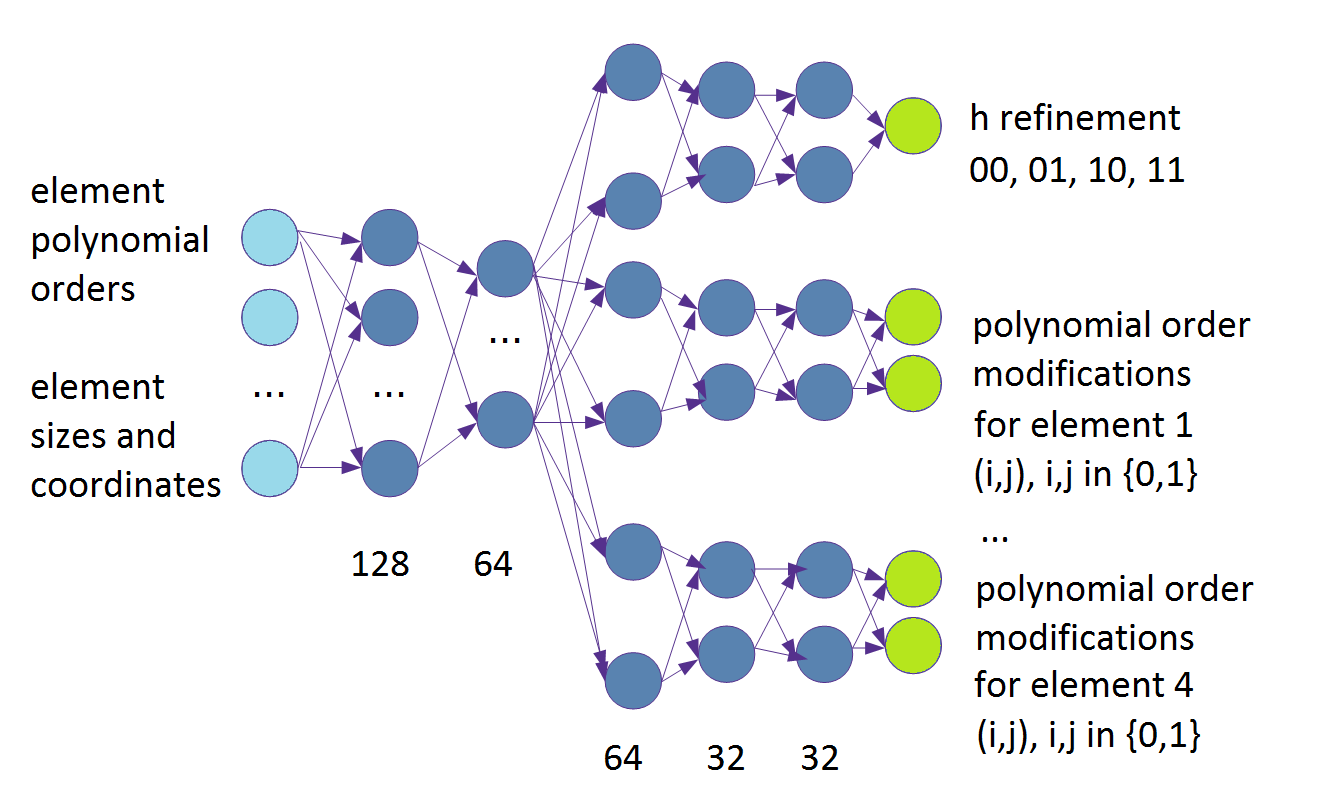}     
         \caption{Deep Neural Network for two dimensional self-adaptive $hp$-FEM. Feed-forward DNN with fully connected layers. ReLU as activation function, softmax as double precision to integer converter (as final activation function in $h$-ref branch). }
\label{fig:DNN2D}
    \end{figure}

    \begin{figure}
        \centering
        \includegraphics[width=0.9\textwidth]{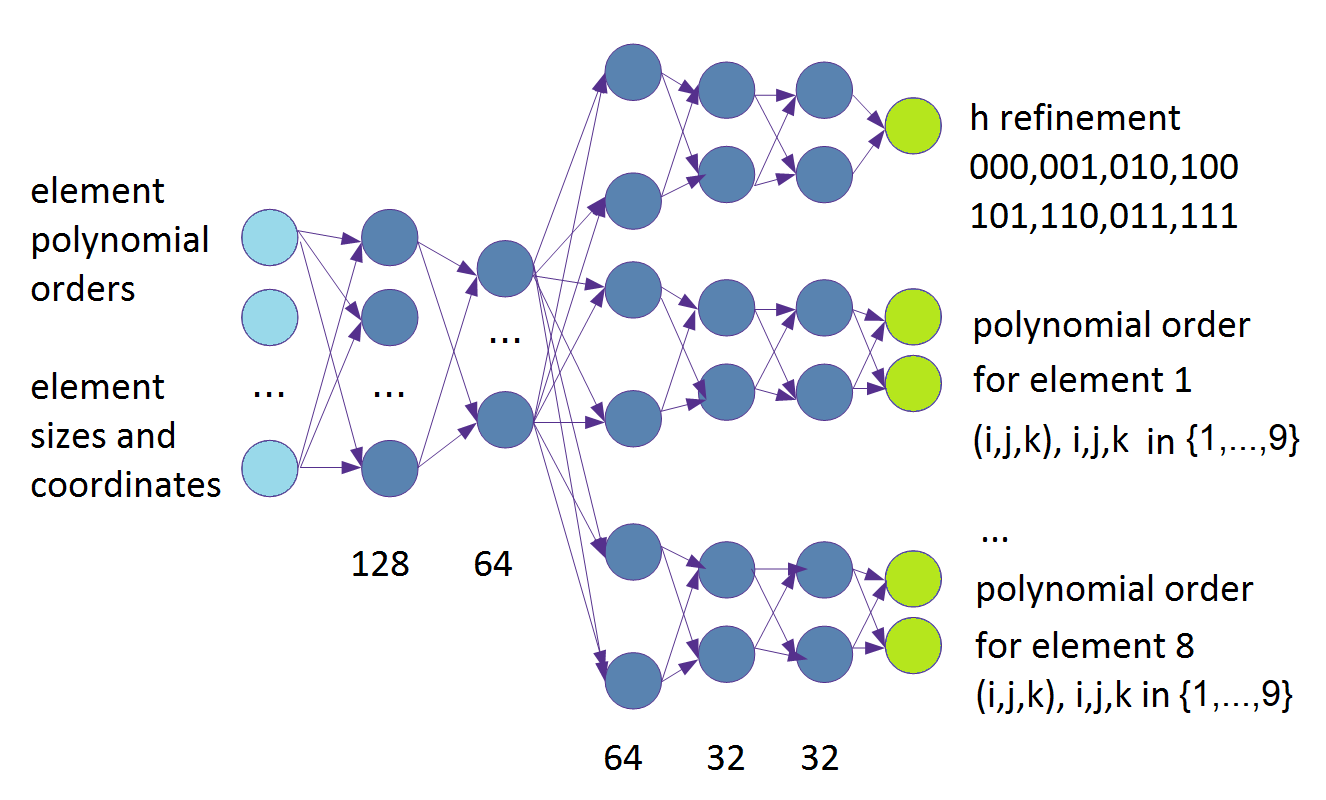}     
         \caption{Deep Neural Network for three dimensional self-adaptive $hp$-FEM. Feed-forward DNN with fully connected layers. ReLU as activation function, softmax as double precision to integer converter (as final activation function in $h$-ref branch). }
\label{fig:DNN3D}
    \end{figure}

We use a feed-forward DNN with fully connected layers \cite{rg1}, presented in Figures \ref{fig:DNN2D}-\ref{fig:DNN3D} The network splits into five branches in 2D and nine branches in 3D, four layers each: the first branch decides about the optimal $h$ refinement, and the remaining branches decide about modifying the polynomial orders - $p$ refinement. Since all possible decisions are encoded as categorical variables, we use cross-entropy as the loss function.

Each layer learns a certain transformation and applies a non-linear activation function – in this case, the ReLU function \cite{rg2}. To prevent overfitting, we use the DropOut mechanism \cite{rg3}, which randomly samples a subset of neurons in each layer and sets their activation to 0 during the training phase. This increases the network’s potential to generalize. The last layer uses linear activation, which is well suited for regression tasks. 

As the output vector dimensionality is distinctly greater than the input data, and additionally highly structural ($h$-adaptation, the orders for eight elements polynomial interpolators), the respective NN cannot have a simple dense MLP architecture. For example, the predicted values of output vector coordinates may be mutually exclusive.  In that case, the volume of the space of irrational solutions is huge and can overwhelm that representing rational ones. Consequently, the NN aproximator will be overfitted, producing mainly useless solutions. Splitting the hidden layers of NN on separate paths, each controlling only those output neurons which represent the same type of information (h-adaptation, element 1 polynomial order, element 2 polynomial order, ..., element 8 polynomial order), limits considerably the volume of solution and parameter spaces. Thus, concentrating the approximator on rational solutions. In the extreme case (all hidden layers are split), we train simultaneously several MLPs on the same data but concentrated on various goals (outputs). Adding more initial not split hidden layers with a considerably larger number of neurons than the input layer allows for finding a better representation of input data.

Now we are ready to introduce DNN-driven $hp$-FEM Algorithm 3 that guides the $hp$ refinements.  It does not employ the two-grid paradigm; its decisions about quasi-optimal $hp$ refinements result from asking the trained DNN.

    \begin{algorithm}[H]
       
        \caption{DNN-driven \textit{hp}-FEM algorithm}
        \KwIn{Initial mesh $T_{init}$, PDE, boundary conditions, \#iterations}
        \KwOut{Optimal mesh, Optimal mesh solution}
        1 iterations = 1 \\
        2 $T_{hp}=T_{init}$ (coarse mesh = initial mesh) \\
        3 \inblue{Solve the coarse mesh problem to obtain $u_{hp}$ (optional)} \\
        4 \If{iterations $==$  \# iterations}{
            \Return{$\left(T_{\frac{h}{2}p+1},u_{hp}\right)$}    
          }
        5 Select optimal refinements $\{ V_{opt}^K\}_{K \in T_{hp}}$ for all elements $K$ by asking DNN \\
$\{ V_{opt}^K = DNN(K)\}_{K \in T_{hp}}$ \ \\
        6 Perform all \textit{hp} refinements from $\{V_{opt}^K\}_{K\in T_{hp}}$ to obtain $T_{opt}$\\
        7 $T_{hp}=T_{opt}$ (coarse mesh = optimal mesh)\\
        8 ++\#iterations \\
        9 \textbf{goto 2}
    \end{algorithm}

Note that the call to the solver in Line 3 is optional, it is only intended to provide the energy norm estimation of the solution $\|u_{hp}\|$ in order to monitor the convergence.

\section{Numerical results}

\label{sec:lshaped}
\begin{figure}[h]
\centering
\begin{subfigure}{0.4\textwidth}
    \includegraphics[width=\textwidth]{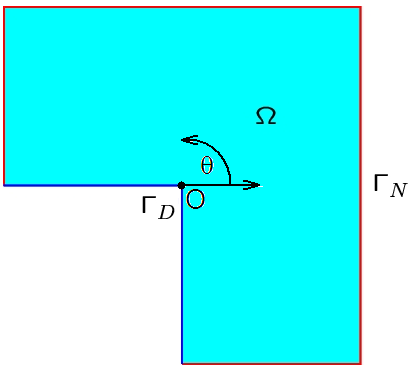}
    \caption{L-shaped}
    \label{fig:Lshape}
\end{subfigure}
\hspace{0.4cm}
\begin{subfigure}{0.48\textwidth}
    \includegraphics[width=\textwidth]{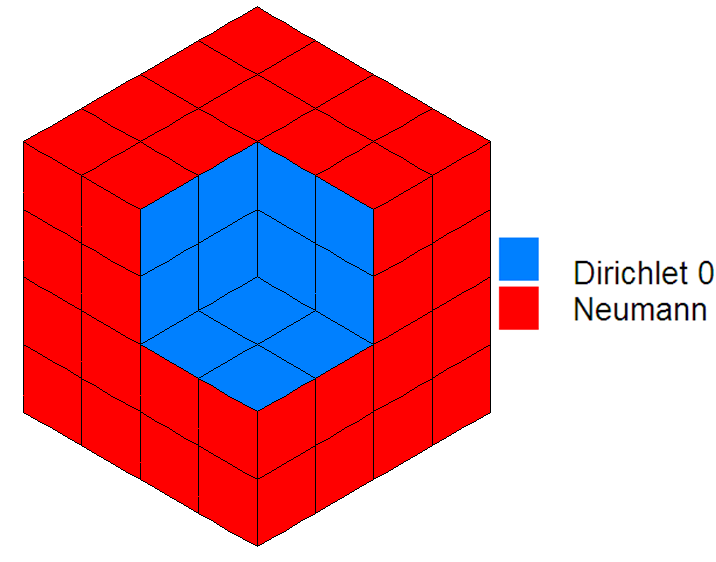}
    \caption{Fichera}
    \label{fig:Fichera}
\end{subfigure}
\caption{Model problem domains}
\label{fig:domains}
\end{figure}

In this section, we present two numerical experiments. 

The first experiment concerns the DNN-driven $hp$-FEM for the two-dimensional L-shaped domain model problem with fixed boundary conditions. 
The first numerical experiment's goal is to verify the convergence of the DNN-driven $hp$-FEM algorithm as compared to the self-adaptive $hp$-FEM algorithm

The second experiment concerns the DNN-driven $hp$-FEM for the 3D Fichera model problem.
The second numerical experiment aims to verify if we can train the DNN using the geometrical data from 12 iterations of the self-adaptive $hp$-FEM to learn the locations of singularities and the distribution patterns of the polynomial orders of approximation. 

\subsection{Verification of the method using two-dimensional L-shaped domain problem}

As the first problem,  we consider the diffusion problem~\eqref{eq:cvf} defined in the L-shaped domain $\Omega = (-1,1)^2\setminus(-1,0)^2$.  We set $\Gamma_D = \{0\} \times (-1,0) \cup (-1,0)\times\{1\}$ and $\Gamma_N  = \text{int}\left(\partial \Omega \setminus \Gamma_N\right)$ as the Dirichlet and Neumann boundaries,  respectively (see Figure~\ref{fig:Lshape}).  As source term,  we consider $g$ such that the analytical solution of Problem~\eqref{eq:cvf} coincides with $u_{\text{ex}}$,  explicitly defined in polar coordinates as:
\begin{equation}
 u_{\text{ex}}(r, \theta) = r^{\frac{2}{3}} \sin\left({\frac{2}{3}}\theta + \frac{\pi}{3}\right). 
\end{equation}
The gradient of the analytical solution is singular at point $(0,0)$ and belongs to the Hilbert space $H^{5/3-\epsilon}(\Omega)$,  for all $\epsilon>0$ (see~\cite{Rivara}).  Last implies in practice that,  to resolve the singularity with high accuracy,  several $h$-refinements and non-uniform distribution of polynomial orders of approximations are required in the layers surrounding the entrant corner (cf.~\cite{Babuska1,Babuska2}.) 

We run the self-adaptive $hp$-FEM algorithm in the model L-shaped domain problem.  
We collect the decisions about the optimal refinements for the coarse mesh elements.
We run 40 iterations of the self-adaptive $hp$-FEM algorithm, and we collect the data summarized in Tables \ref{tab:input2D}-\ref{tab:output2D} as the dataset for training.
We collect around 10,000 samples for the dataset by recording the decisions about the optimal $hp$ refinements performed for particular finite elements.

Having the DNN trained, we continue with the L-shaped domain problem computations. We run 10 iterations of the DNN-driven $hp$-FEM.
In this 2D problem, we can also continue with the self-adaptive $hp$-FEM algorithm (running without the DNN). Thus, we compare last 10 iterations of our DNN-driven $hp$-FEM against the last 10 iterations of the self-adaptive $hp$-FEM algorithm. 

\begin{figure}
    \centering
    
    \includegraphics[width=0.9\textwidth]{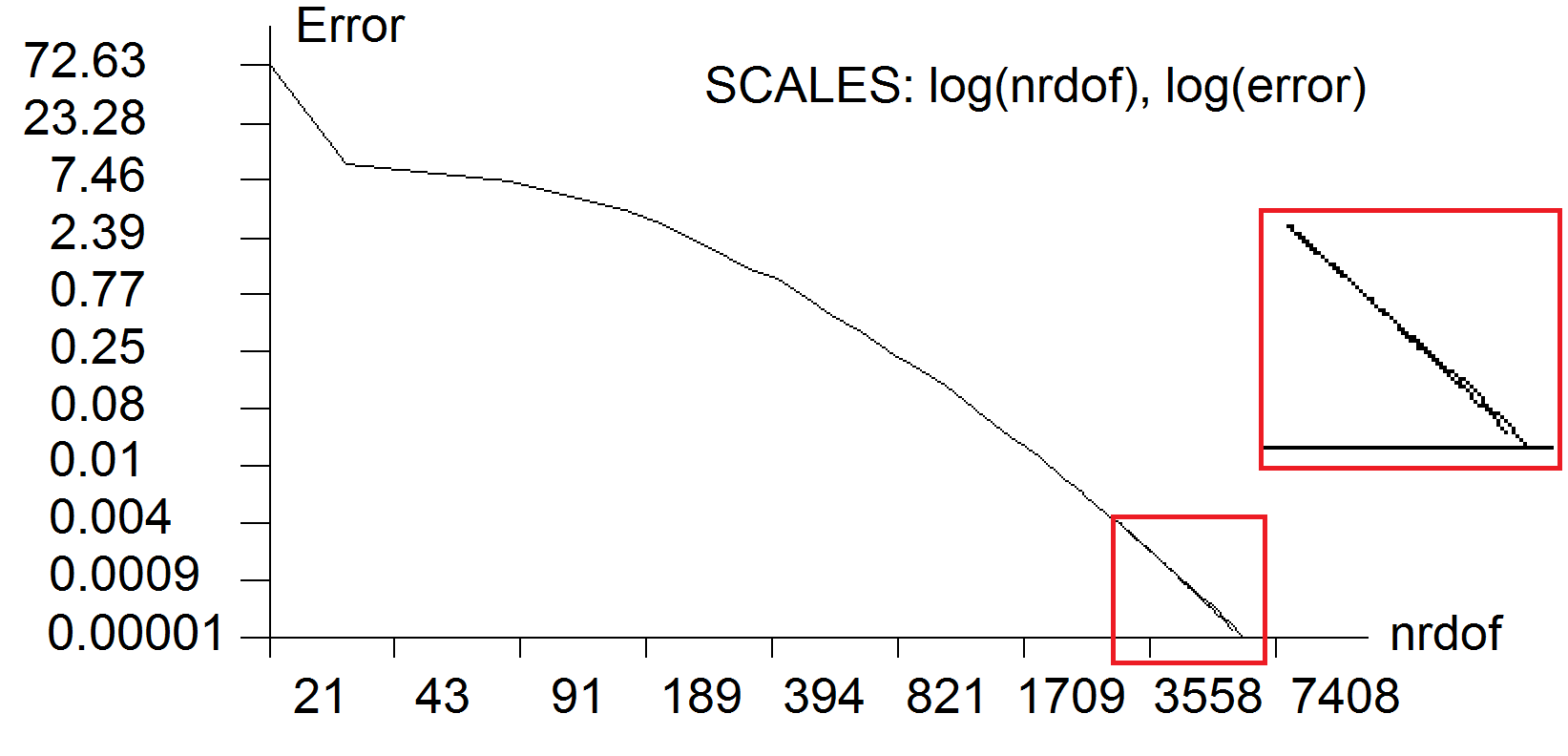}

    \caption{The comparison of 50 iterations of self-adaptive $hp$-FEM and hybrid 40 iterations of self-adaptive and 10 iterations of DNN-driven $hp$-FEM on original L-shaped domain. }
        \label{fig34}

\end{figure}

\begin{figure}
    \centering
    \includegraphics[width=0.49\textwidth]{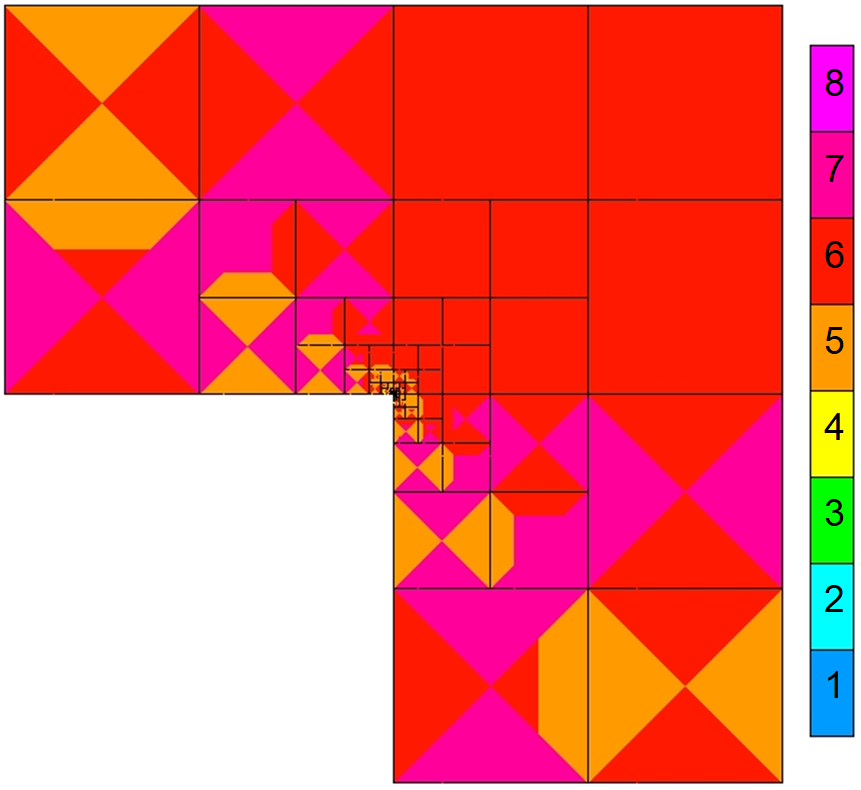}
    \includegraphics[width=0.45\textwidth]{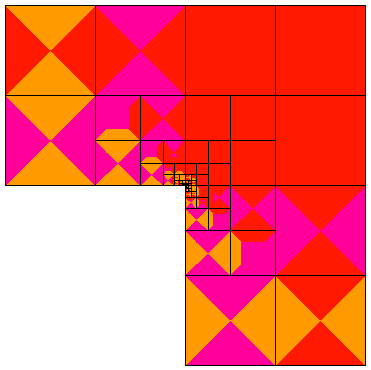}\\
    
    \caption{The mesh provided by 50 iterations of the self-adaptive $hp$-FEM algorithm (left panel) \\ and by 40 iterations of self-adaptive $hp$-FEM followed by 10 iterations of the deep learning-driven $hp$-FEM algorithm (right panel).}
    
    \label{fig:comparison1}
\end{figure}

\begin{figure}
    \centering
    \includegraphics[width=0.49\textwidth]{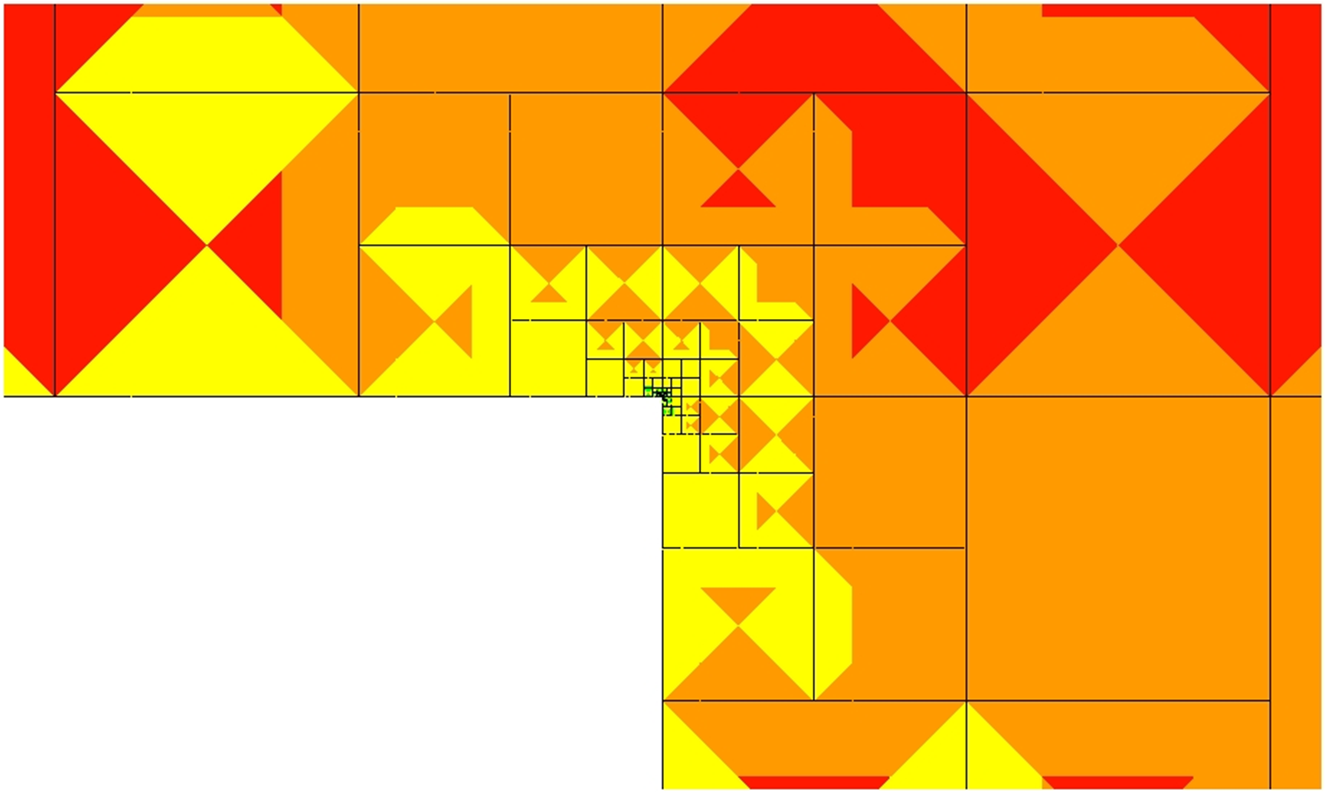}
    \includegraphics[width=0.49\textwidth]{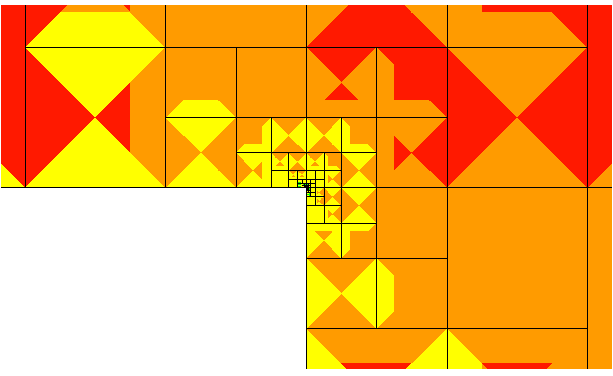}
 \caption{The mesh provided by 50 iterations of the self-adaptive $hp$-FEM algorithm (left panel) \\ and by 40 iterations of self-adaptive $hp$-FEM followed by 10 iterations of the deep learning-driven $hp$-FEM algorithm (right panel).  Zoom 100 X}
    \label{fig:comparison2}
\end{figure}

\begin{figure}
    \centering
    \includegraphics[width=0.49\textwidth]{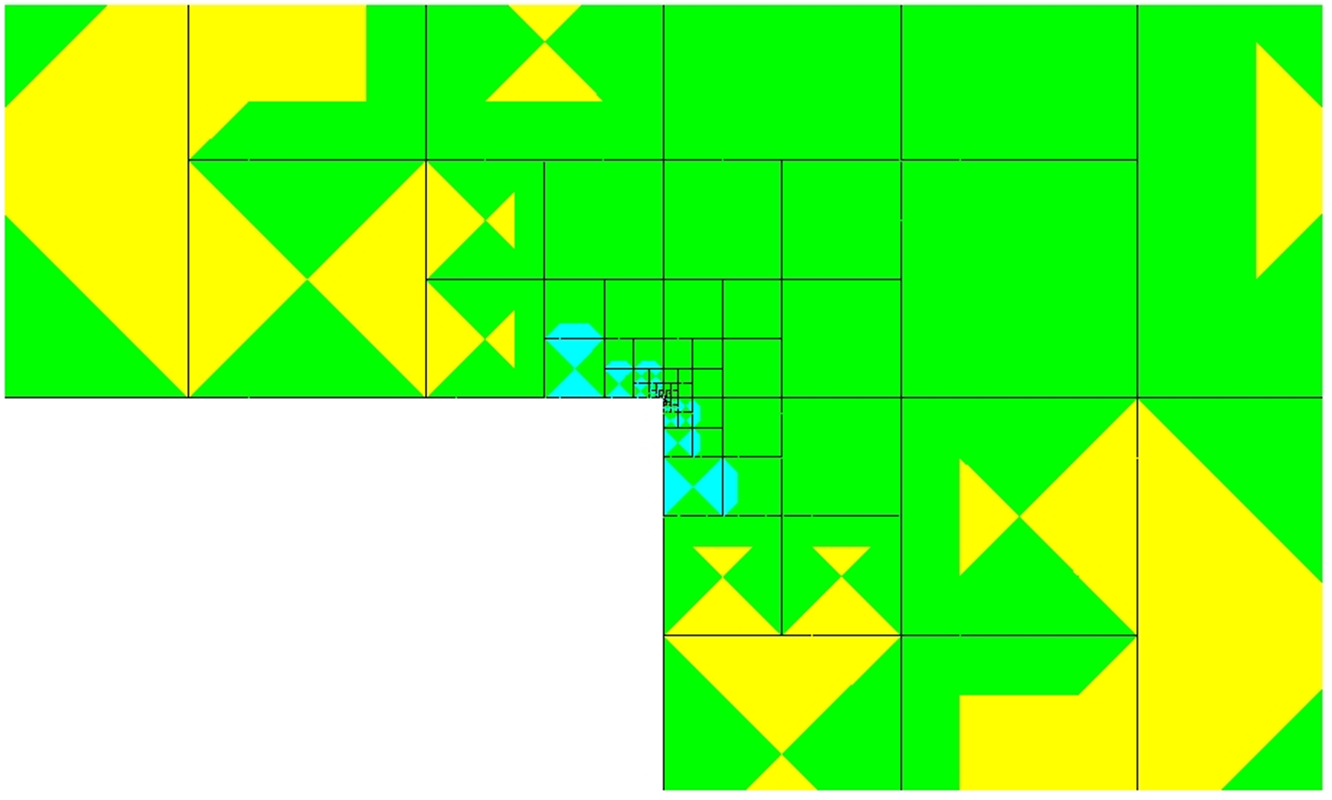}
   \includegraphics[width=0.49\textwidth]{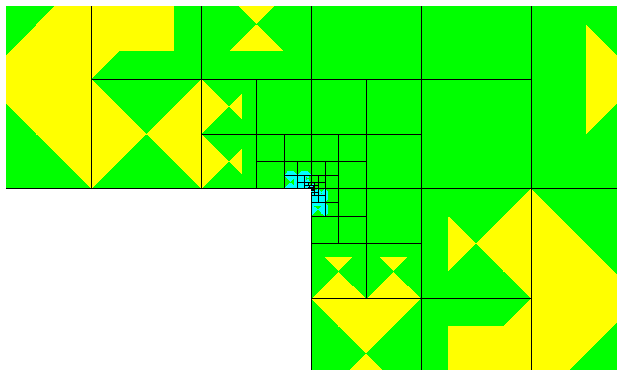}    
\caption{The mesh provided by 50 iterations of the self-adaptive $hp$-FEM algorithm (left panel) \\ and by 40 iterations of self-adaptive $hp$-FEM followed by 10 iterations of the deep learning-driven $hp$-FEM algorithm (right panel).
 Zoom 10,000 X }   
    \label{fig:comparison3}
\end{figure}

\begin{figure}
    \centering
    \includegraphics[width=0.49\textwidth]{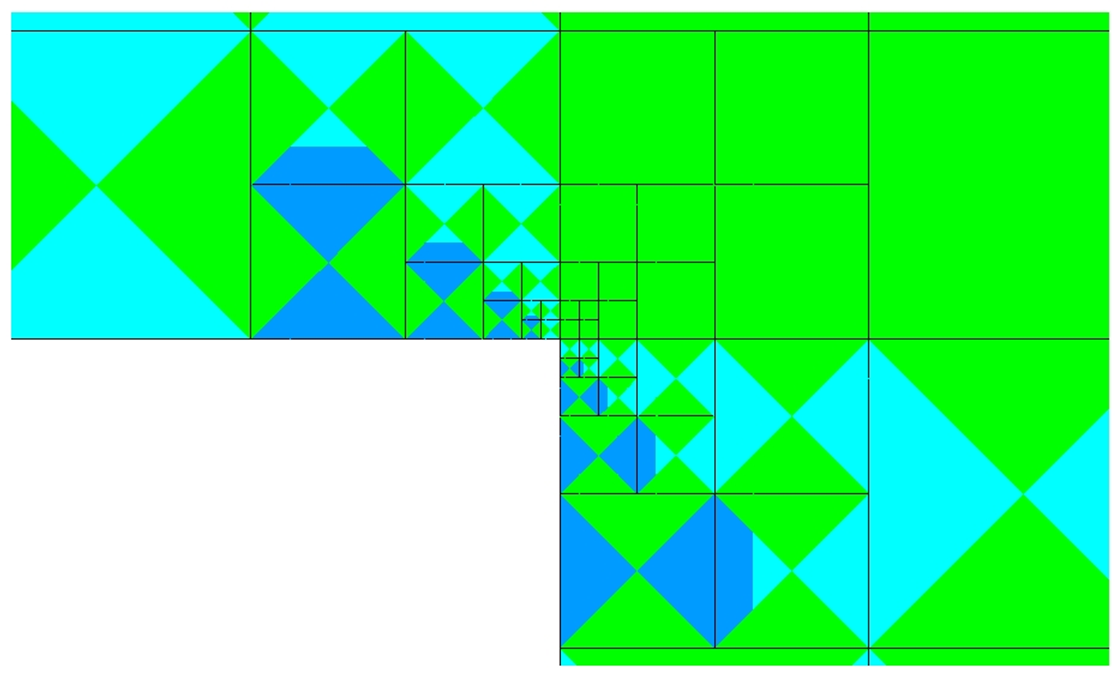}
    \includegraphics[width=0.49\textwidth]{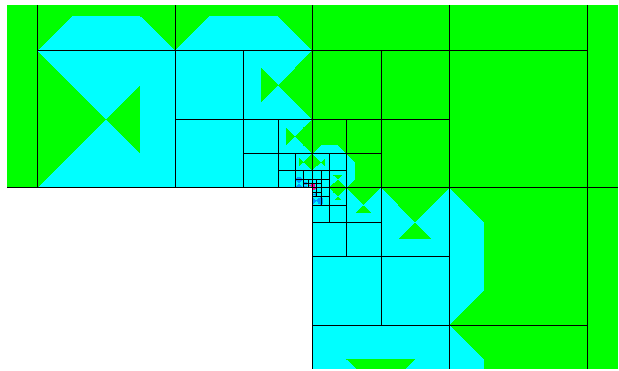}    
\caption{The mesh provided by 50 iterations of the self-adaptive $hp$-FEM algorithm (left panel) \\ and by 40 iterations of self-adaptive $hp$-FEM followed by 10 iterations of the deep learning-driven $hp$-FEM algorithm (right panel).  \\ Zoom 1,000,000 X }
    \label{fig:comparison4}
\end{figure}

We compare the exponential convergence of the self-adaptive $hp$-FEM and the DNN-driven $hp$-FEM algorithms in Figure \ref{fig34}. We can see that both methods deliver similar exponential convergence, and we have trained the DNN very well.
We compare in Figures \ref{fig:comparison1}-\ref{fig:comparison4} the final mesh generated by both methods.
The last 10 iterations of the DNN-driven $hp$-FEM algorithm performs refinements similar to the last 10 iterations of the self-adaptive $hp$-FEM. The difference concerns very small elements, with a diameter of $10^{-6}$, depicted in Figure \ref{fig:comparison4}.
We conclude that in two dimensions, we can effectively train the DNN, so it predicts the optimal refinements and delivers exponential convergence.

\begin{table}
\centering
\begin{tabular}{cc}
{\bf Input parameter} & {\bf Data dimensionality } \\
 \hline
Element coordinates  & 4 \\
$(x_1,y_1),(x_2,y_2)$ & double precision\\
\hline
Element dimensions & 2 \\
$(d_x,d_y)$ & double precision\\
\hline
polynomial order of approximation for element & 2 \\
$(p_x,p_y)$ where $p_x,p_y \in \{1,...,9\} $ & integer \\
\hline
\vspace{5pt}
\end{tabular}
\caption{The input for the DNN concerns the local element data from the coarse mesh and one global value - the maximum over all the elements of the norms of the solution. The total size of the input data is 8 double precision values.}
\label{tab:input2D}
\end{table}

\begin{table}
\centering
\begin{tabular}{cc}
{\bf Output parameter} & {\bf Data dimensionality } \\
 \hline
$h$ refinement flag & 3 \\
$(href_x,href_y)$, where $href_x,href_y \in \{0,1\}$  & boolean \\
\hline
polynomial order of approximation for element 1 & 2 \\
$(p_x^1,p_y^1)$ where $p_x^1,p_y^1 \in \{1,...,9\} $ & integer \\
\hline
polynomial order of approximation for element 2 & 2 \\
$(p_x^2,p_y^2)$ where $p_x^2,p_y^2 \in \{1,...,9\} $ & integer \\
\hline
polynomial order of approximation for element 3 & 2 \\
$(p_x^3,p_y^3)$ where $p_x^2,p_y^2 \in \{1,...,9\} $ & integer \\
\hline
polynomial order of approximation for element 4 & 2 \\
$(p_x^4,p_y^4)$ where $p_x^2,p_y^2 \in \{1,...,9\} $ & integer \\
\hline
\vspace{5pt}
\end{tabular}
\caption{The output from the DNN concerns the $h$-refinement flags, and the polynomial orders of approximation for up to 4 son elements. For the case when the refinement is not needed, ($href_x href_y=00)$ the element orders are encoded in $(p_x^1,p_y^1)$ and the other $(p_x^i,p_y^i)=(0,0)$ for $i=2,3,4$. For the case when the horizontal or vertical refinement is selected, ($href_x href_y \in \{10,01\}$ the two son element orders are encoded in $(p_x^1,p_y^1),(p_x^2,p_y^2)$ and the other $(p_x^i,p_y^i)=(0,0)$ for $i=3,4$. }
\label{tab:output2D}
\end{table}


\subsection{Three-dimensional Fichera model problem}

This 3D problem was introduced by Gaetano Fichera \cite{Fichera} as the benchmark to test the convergence of 3D finite element method codes. It has a point singularity in the central point (0,0,0) and three edge-singularities around the three central edges.

This problem is the generalization of the L-shaped domain model problem into three dimensions.
We seek the temperature scalar field over
$\Omega = (-1,1)^3\setminus(0,1)^3$.  We set $\Gamma_D =(0,1)\times \{ 0\}\times(0,1) \cup \{0\}\times (-1,0)\times(0,1) \cup (0,1)\times (-1,0) \times \{0\} $ and $\Gamma_N  = \text{int}\left(\partial \Omega \setminus \Gamma_N\right)$ as the Dirichlet and Neumann boundaries,  respectively (see Figure~\ref{fig:Fichera}). 
The Neumann boundary condition is obtained from the superposition of three solutions of the two-dimensional L-shaped domain model problem. The details of the formulation are provided in \cite{Rachowicz}.
\begin{equation}
\Delta u = 0 \text{ in } \Omega \qquad 
u = 0 \text{ on } \Gamma_D \qquad 
\frac{\partial u}{\partial n} =g \text{ on } \Gamma_N  
\end{equation}
We consider the weak formulation
\begin{equation}
\label{eq:Fichera}
b(u,v) =l(v) \forall v\in V \quad b(u,v) = \int_{\Omega}\nabla u \nabla v dx \quad l(v) = \int_{\Gamma_N}gv dS  
\end{equation}
and introduce proper Sobolev space
\begin{equation}
V=\Bigg\{v \in L^2(\Omega): \int_\Omega ||v||^2 + ||\nabla v||^2 < \infty: tr(v) = 0 \text{ on } \Gamma_D\Bigg\}
\end{equation}
To resolve the one point and three edge singularities with high accuracy, this problem requires several $h$ refinements towards the edges and non-uniform distribution of polynomial orders of approximations in the layers surrounding the central point and the three edges.

\begin{figure}[h]
        \centering
        \includegraphics[width=0.2\textwidth]{coarse.png}     $\rightarrow$   \includegraphics[width=0.2\textwidth]{second.png}     $\rightarrow$   \includegraphics[width=0.2\textwidth]{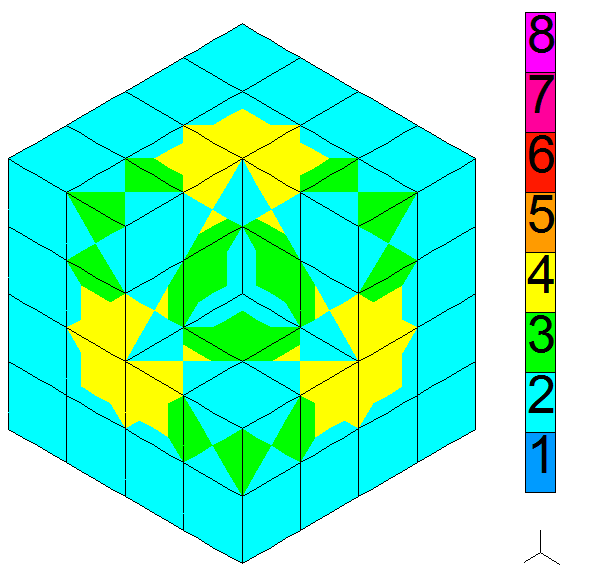}     $\rightarrow$   \includegraphics[width=0.2\textwidth]{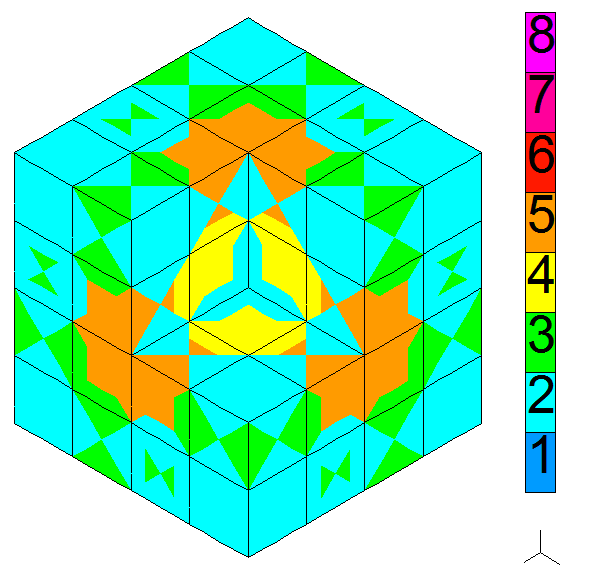}     $\rightarrow$\\ \includegraphics[width=0.2\textwidth]{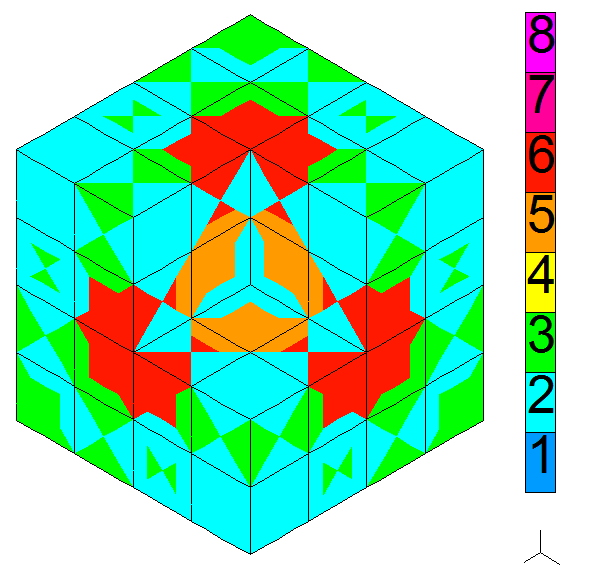}     $\rightarrow$   \includegraphics[width=0.2\textwidth]{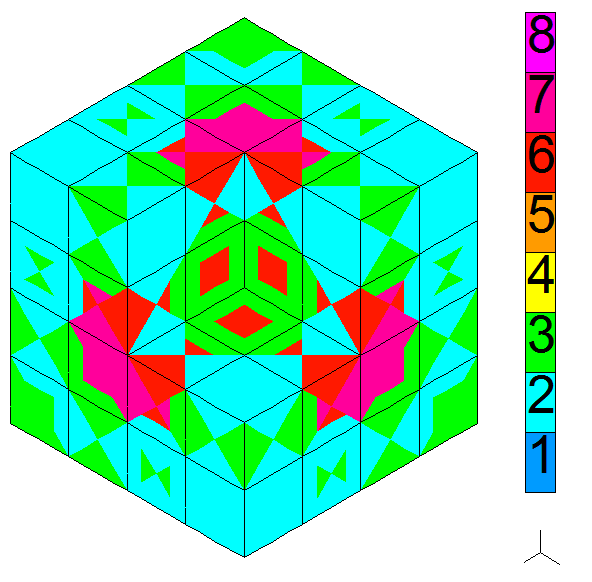}     $\rightarrow$   \includegraphics[width=0.2\textwidth]{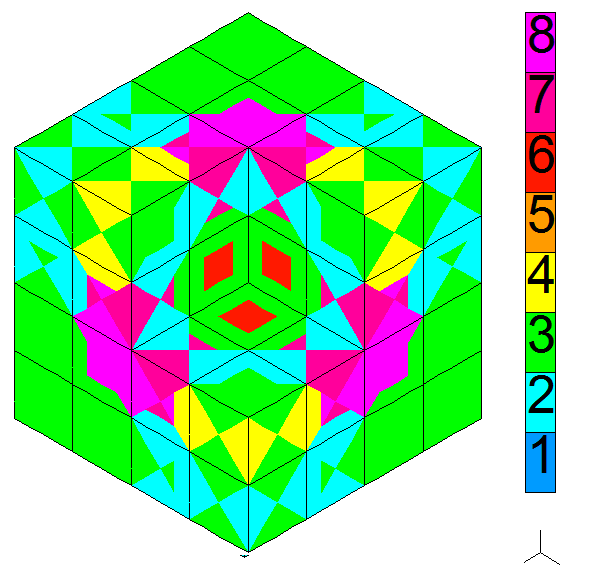}     $\rightarrow$   \includegraphics[width=0.2\textwidth]{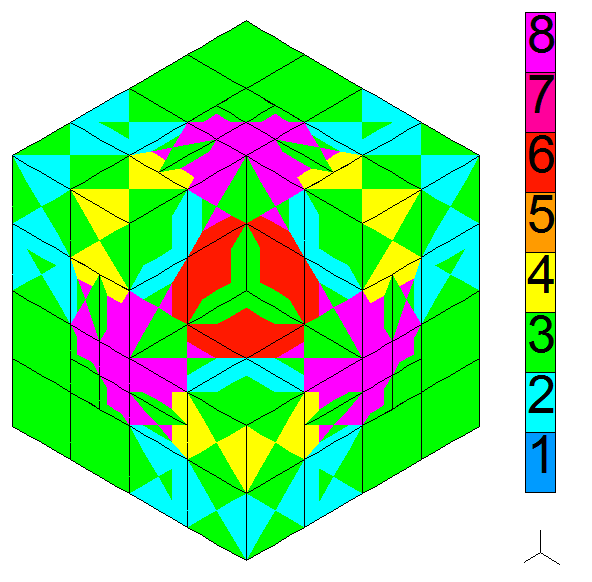}     $\rightarrow$\\
\includegraphics[width=0.2\textwidth]{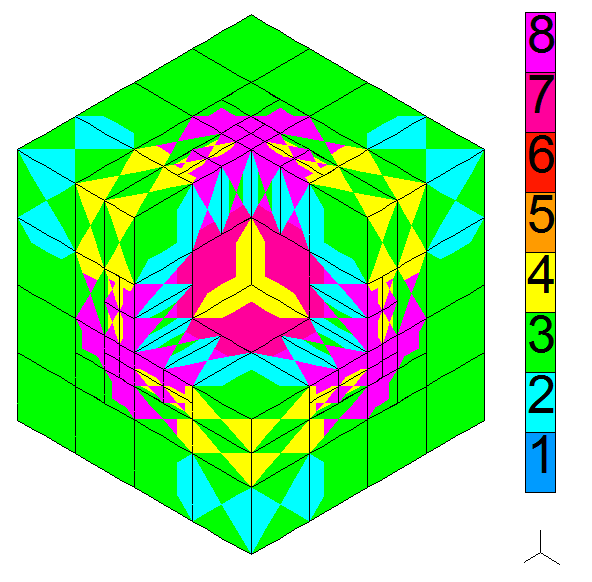}     $\rightarrow$   \includegraphics[width=0.2\textwidth]{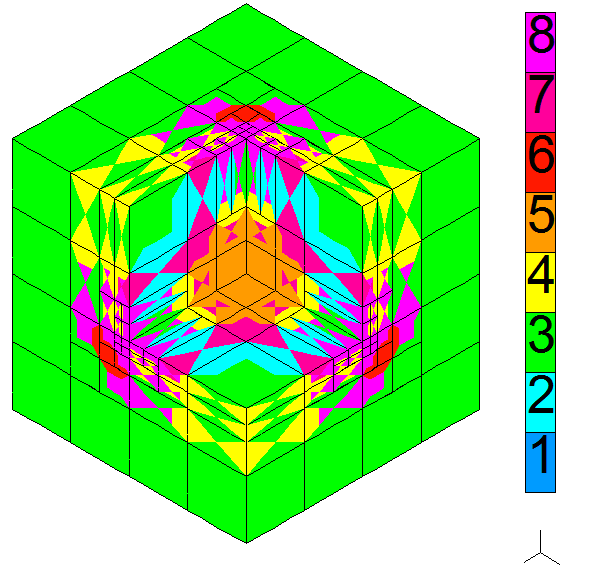}     $\rightarrow$   \includegraphics[width=0.2\textwidth]{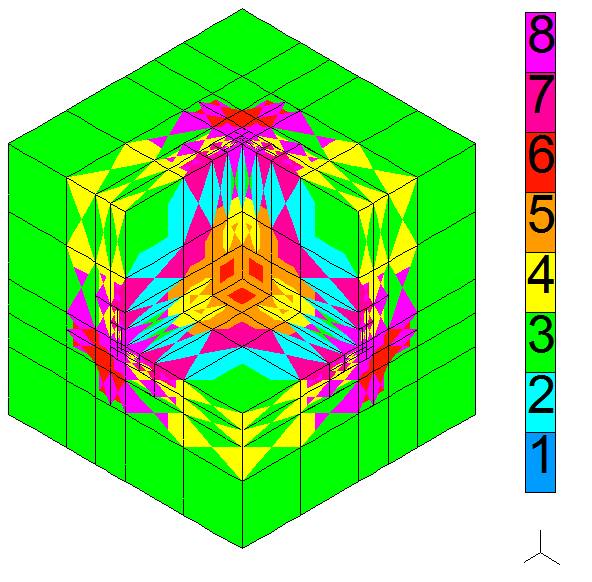} $\rightarrow$     
\includegraphics[width=0.2\textwidth]{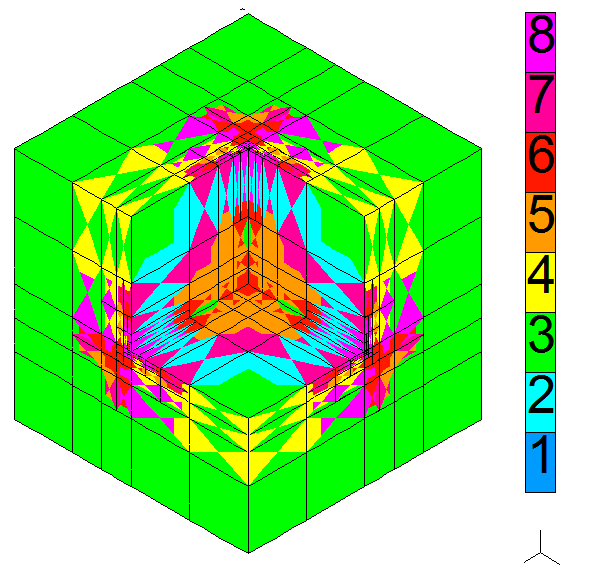} 
        \caption{Sequence of $hp$ adaptive meshes generated by the self-adaptive $hp$ finite element method algorithm}        
\label{fig:Ficheragrids}
    \end{figure}

In this numerical experiment, we run Algorithms 1 and 2 of the self-adaptive $hp$-FEM algorithm for 12 iterations of the Fichera model problem. Figure \ref{fig:Ficheragrids} presents the sequence of generated meshes. 
We collect around 2000 samples for the dataset.
The input and the output employed for the training of the DNN are summarized in Tables \ref{tab:input3D}-\ref{tab:output3D}.
The DNN learns the location of the singularities, where the $h$ refinements are needed, based on the geometrical coordinates and dimensions of the finite elements. It also learns the $p$ refinement patterns, using the geometrical data and the element distribution of the polynomial approximation orders.

\begin{table}
\centering
\begin{tabular}{cc}
{\bf Input parameter} & {\bf Data dimensionality } \\
 \hline
Element coordinates  & 6 \\
$(x_1,y_1,z_1),(x_2,y_2,z_2)$ & double precision\\
\hline
Element dimensions & 3 \\
$(d_x,d_y,d_z)$ & double precision\\
\hline
Polynomial orders of approximation for element interior & 3 \\
$(p_x,p_y,p_z)$ & integer\\
\hline
\vspace{5pt}
\end{tabular}
\caption{The input for the DNN concerns the geometrical data from the coarse mesh and one local element distribution of polynomial orders. The total size of the input data is 9 double precision and 3 integer values.}
\label{tab:input3D}
\end{table}

\begin{table}
\centering
\begin{tabular}{cc}
{\bf Output parameter} & {\bf Data dimensionality } \\
 \hline
$h$ refinement flag & 3 \\
$(href_x,href_y,href_z)$, where $href_x,href_y,href_z \in \{0,1\}$  & boolean \\
\hline
polynomial order of approximation for element 1 & 2 \\
$(p_x^1,p_y^1)$ where $p_x^1,p_y^1 \in \{1,...,9\} $ & integer \\
\hline
polynomial order of approximation for element 2 & 2 \\
$(p_x^2,p_y^2)$ where $p_x^2,p_y^2 \in \{1,...,9\} $ & integer \\
\hline
polynomial order of approximation for element 3 & 2 \\
$(p_x^3,p_y^3)$ where $p_x^2,p_y^2 \in \{1,...,9\} $ & integer \\
\hline
polynomial order of approximation for element 4 & 2 \\
$(p_x^4,p_y^4)$ where $p_x^2,p_y^2 \in \{1,...,9\} $ & integer \\
\hline
polynomial order of approximation for element 5 & 2 \\
$(p_x^4,p_y^4)$ where $p_x^2,p_y^2 \in \{1,...,9\} $ & integer \\
\hline
polynomial order of approximation for element 6 & 2 \\
$(p_x^4,p_y^4)$ where $p_x^2,p_y^2 \in \{1,...,9\} $ & integer \\
\hline
polynomial order of approximation for element 7 & 2 \\
$(p_x^4,p_y^4)$ where $p_x^2,p_y^2 \in \{1,...,9\} $ & integer \\
\hline
polynomial order of approximation for element 8 & 2 \\
$(p_x^4,p_y^4)$ where $p_x^2,p_y^2 \in \{1,...,9\} $ & integer \\
\hline
\vspace{5pt}
\end{tabular}
\caption{The output from the DNN concerns the $h$-refinement flags, and the polynomial orders of approximation for up to 8 son elements. For the case when the refinement is not needed, $(href_x= href_y=href_z=0)$ the element orders are encoded in $(p_x^1,p_y^1,p_z^1)$ and the other $(p_x^i,p_y^i,p_z^i)=(0,0,0)$ for $i=2,...,8$. For the case of one-directional refinements $(href_x href_y href_z\in \{100,010,001\}$ the two son element orders are encoded in $(p_x^1,p_y^1,p_z^1),(p_x^2,p_y^2,p_z^2)$ and the other $(p_x^i,p_y^i,p_z^i)=(0,0,0)$ for $i=3,...,8$. For the case of two-directional refinements $(href_x href_y href_z\in \{110,101,011\}$ the four son element orders are encoded in $(p_x^i,p_y^i,p_z^i)$ where $i=1,2,3,4$, and the other $(p_x^i,p_y^i,p_z^i)=(0,0,0)$ for $i=5,...,8$.}
\label{tab:output3D}
\end{table}

\begin{figure}
        \centering
\includegraphics[width=0.9\textwidth]{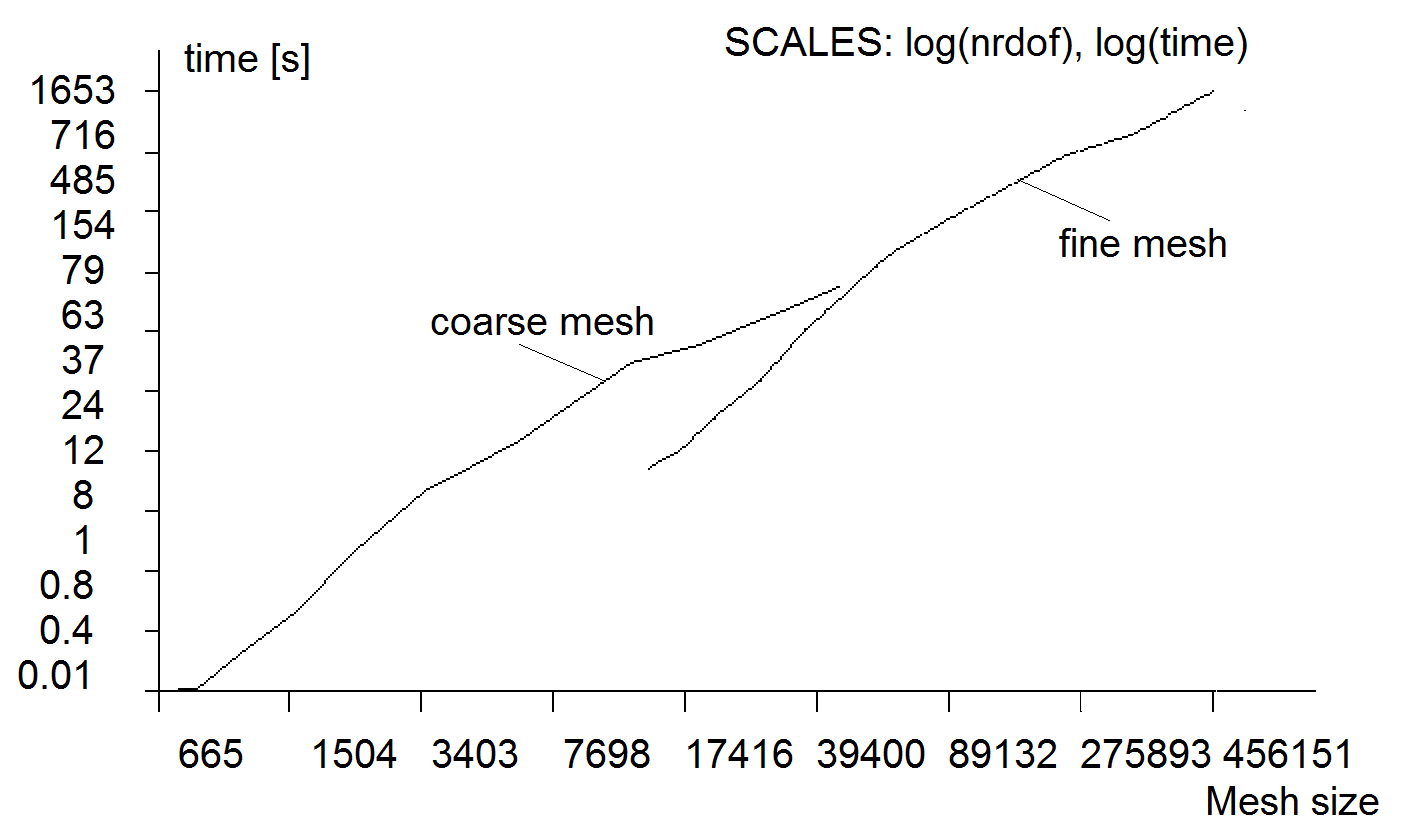}     

         \caption{The computational cost of the training process. Execution time of twelve iterations on a laptop with 8GB of RAM and 2.5 GHz processor.}
        \label{fig:time}
        
    \end{figure}

\begin{figure}
        \centering
\includegraphics[width=0.9\textwidth]{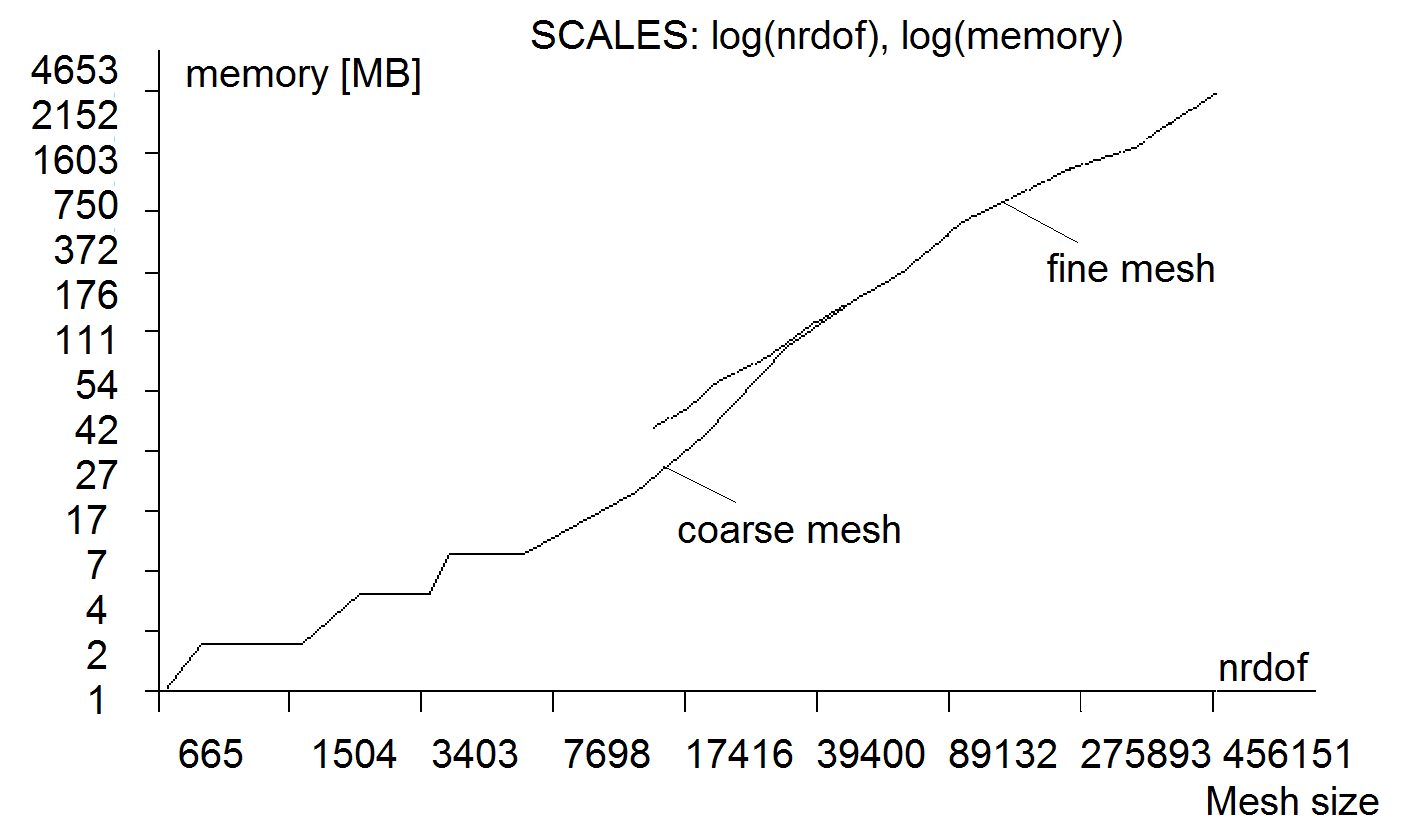}     

         \caption{The memory cost of the training process. Memory consumption of the solver after twelve iterations on a laptop with 8GB of RAM and 2.5 GHz processor.}
        \label{fig:memory}
    \end{figure}

The Algorithm 1 of the self-adaptive $hp$-FEM algorithm employs the two grids paradigm. It solves the problem over the coarse and the fine mesh, using the MUMPS solver \cite{MUMPS1,MUMPS2,MUMPS3}.
In this 3D model problem, the memory requirements and the execution times for the coarse and fine mesh problems are summarized in Figures \ref{fig:time}-\ref{fig:memory}.

\begin{figure}
        \centering
        \includegraphics[width=0.2\textwidth]{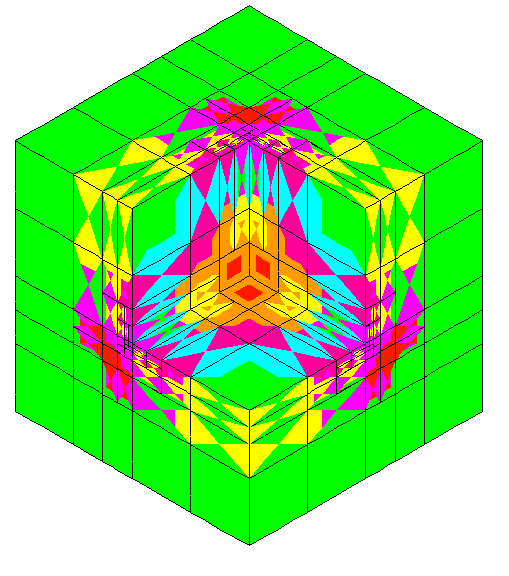}     $\rightarrow$   \includegraphics[width=0.2\textwidth]{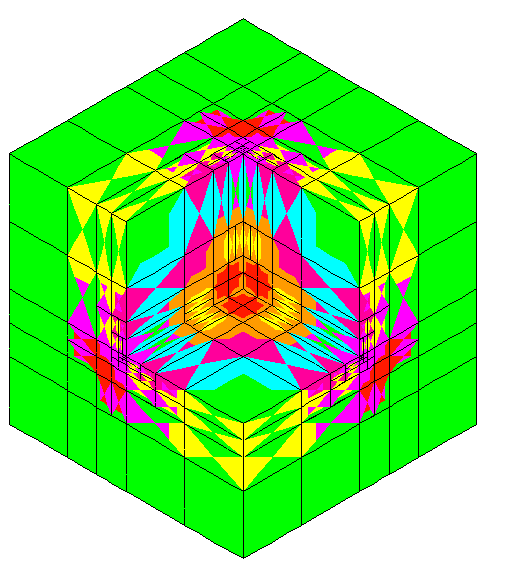}     $\rightarrow$   \includegraphics[width=0.2\textwidth]{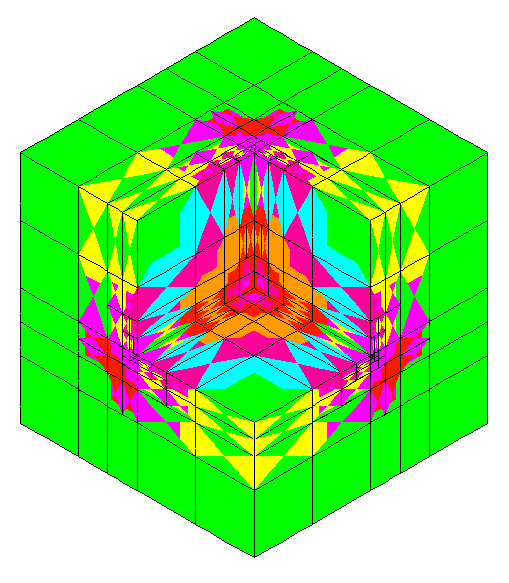}     $\rightarrow$   \includegraphics[width=0.2\textwidth]{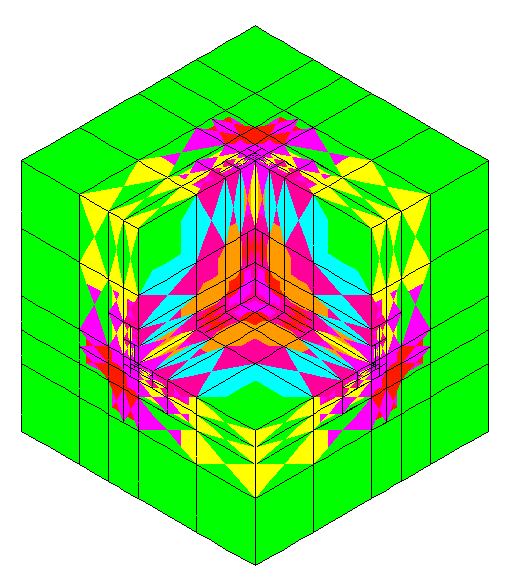}     $\rightarrow$   \includegraphics[width=0.2\textwidth]{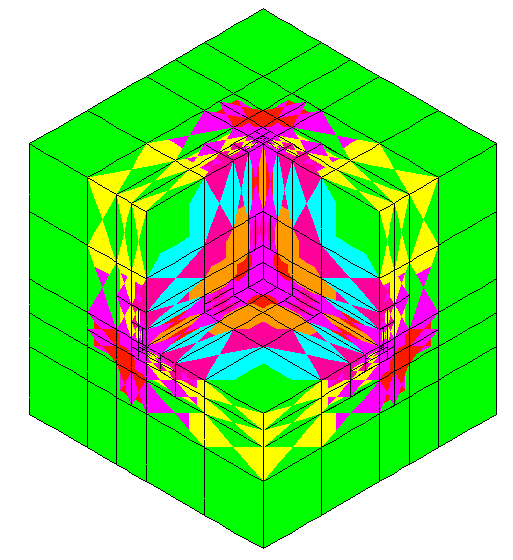}     $\rightarrow$           
\includegraphics[width=0.2\textwidth]{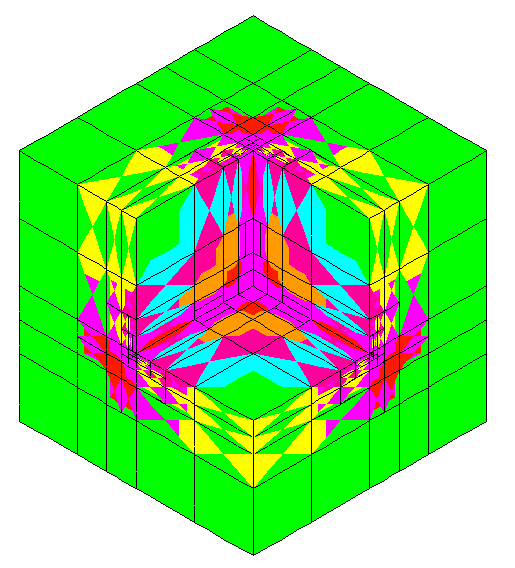}    
  \caption{Sequence of $hp$ adaptive meshes generated by DNN-driven algorithm after 12 iterations of $hp$-adaptive FEM}
\label{fig:DNNFichera}
   \end{figure}

We stop the self-adaptive $hp$-FEM algorithm after 12 iterations. We use the collected 2000 samples to train the DNN. This time no calling the solver, we continue the refinements by using Algorithm 3 of the DNN-driven $hp$-FEM. The further sequence of meshes generated by the DNN-driven $hp$-FEM algorithm is presented in Figure \ref{fig:DNNFichera}.

\begin{figure}
        \centering
        \includegraphics[width=0.8\textwidth]{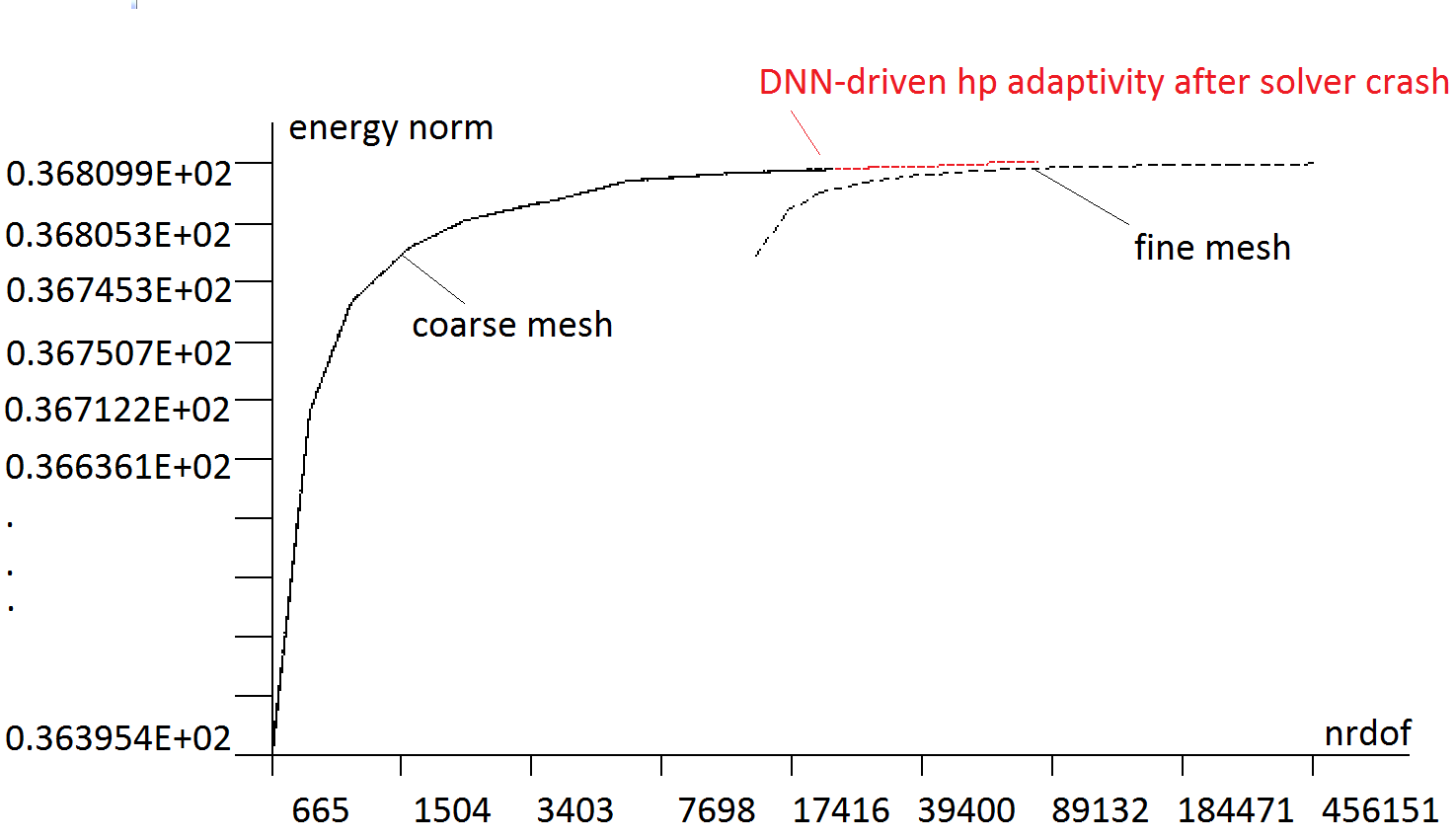}     
          \caption{Convergence of energy of the solution for self-adaptive $hp$-FEM on coarse mesh, fine mesh, and the DNN-driven $hp$-FEM extension}
        \label{fig:energy}
    \end{figure}

\begin{figure}
        \centering
        \includegraphics[width=1.0\textwidth]{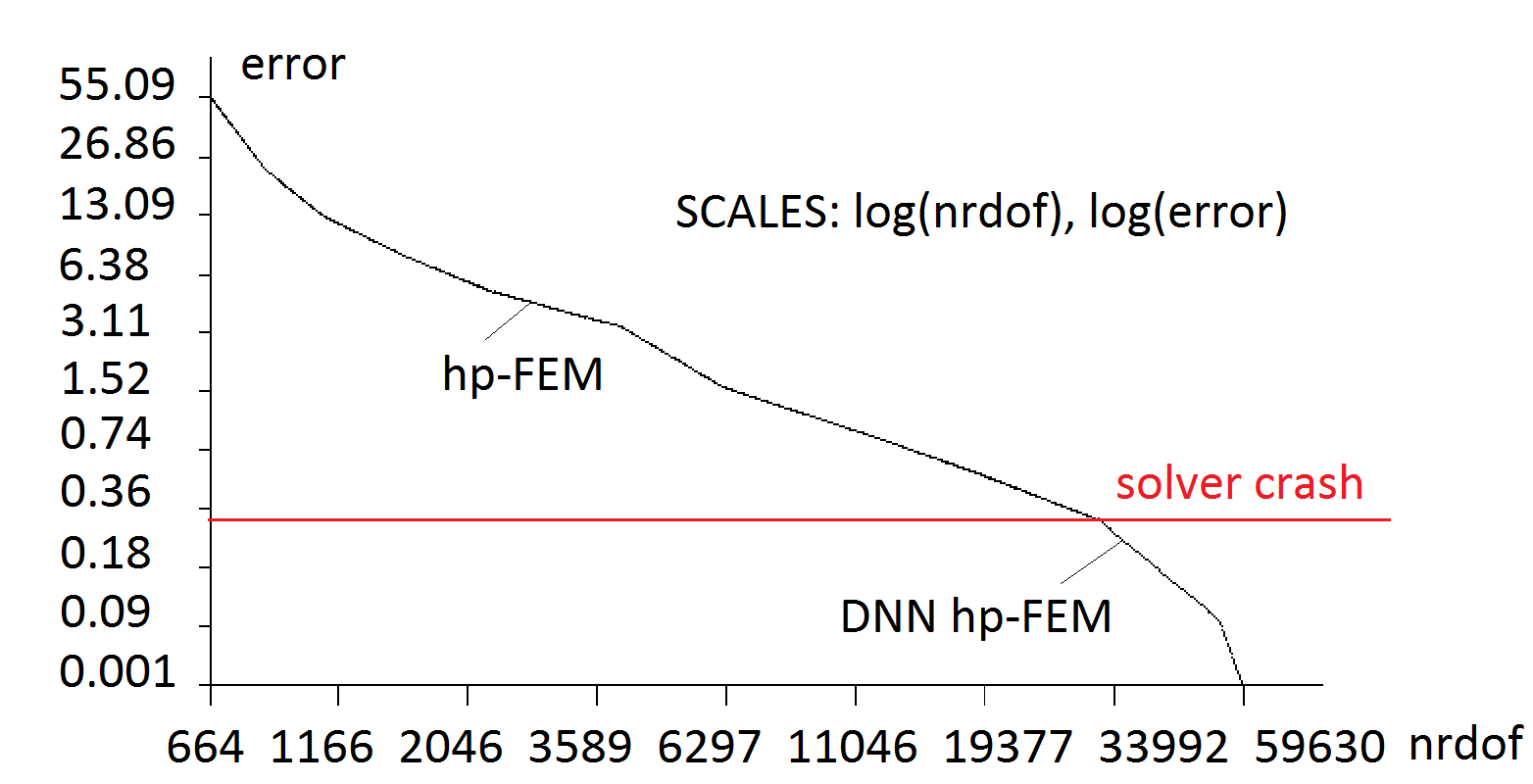}     
          \caption{Convergence of relative error $100\times \frac{|energy - energy\_fine|}{|energy\_fine|}$ computed with respect to the energy of the finest mesh solution for self-adaptive $hp$-FEM, and the continuation with DNN-driven $hp$-FEM extension}
        \label{fig:relative}
    \end{figure}

In order to verify the convergence of the DNN-driven $hp$-FEM, we estimate the energy norm of the solution on the sequence of coarse grids generated by either self-adaptive $hp$-FEM or DNN-driven $hp$-FEM, and the fine grids for the case of th self-adaptive $hp$-FEM. The convergence plots are summarized in Figure \ref{fig:energy}. By using the energy norm of the finest solution obtained on the fine mesh from twelve iterations of the self-adaptive $hp$-FEM algorithm, we can plot the convergence of the relative error of the self-adaptive $hp$-FEM and the DNN-driven $hp$-FEM, presented in Figure \ref{fig:relative}. The red line denotes the moment when we finish the self-adaptive $hp$-FEM (e.g., the fine mesh solver runs out of resources), and we continue with the DNN-driven $hp$-FEM. As a result, we obtain a very nice exponential convergence. To understand these results, we must remember that we compute the relative error with respect to the energy norm of the finest mesh solution.

\subsection{Deep Neural Network Training}

We train the network for 200 epochs using the Adam optimizer \cite{rg4}, which minimizes the Mean Squared Error cost function. The error backpropagation algorithm uses information from the optimizer to correct the network’s weights. We can modify the magnitude of that correction by changing the learning rate constant. In addition, we regularize the network by reducing the learning rate whenever the loss function minimization slows down. We also use early stopping to stop the training procedure before completing all 200 training epochs if the difference in results between epochs is insignificant and cannot be improved by the learning rate reduction scheme.
We used the TensorFlow 1.13 and CUDA 9 library to train the network on an nVidia Tesla v100 GPU with 32GB VRAM and 64GB general-purpose RAM.

\begin{figure}
    \centering
    \includegraphics[width=0.49\textwidth]{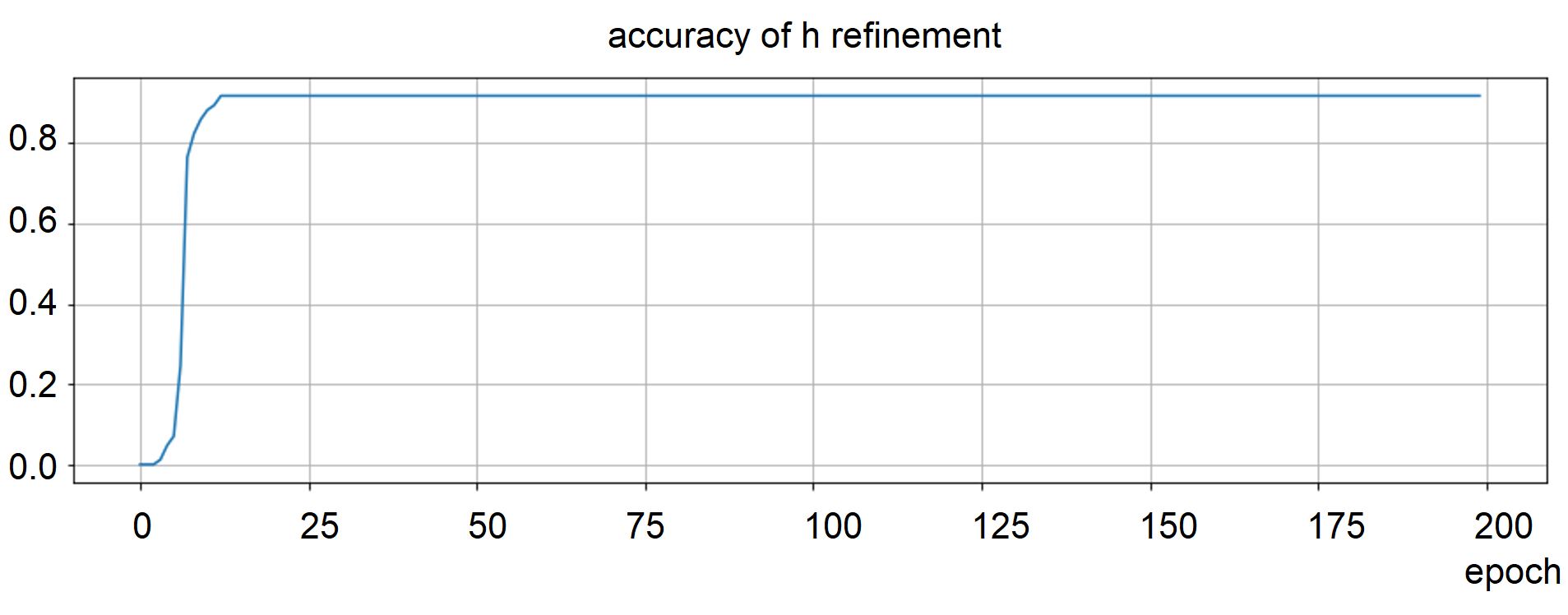} \\
    \includegraphics[width=0.49\textwidth]{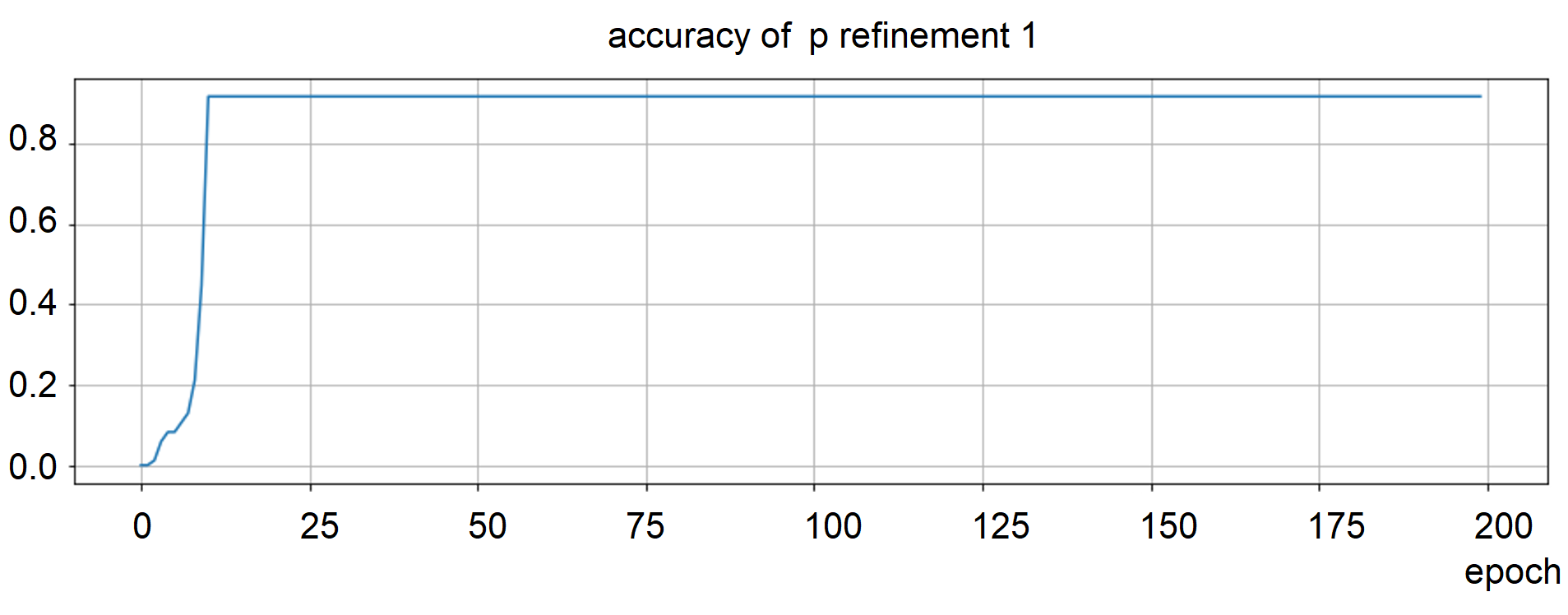}    \includegraphics[width=0.49\textwidth]{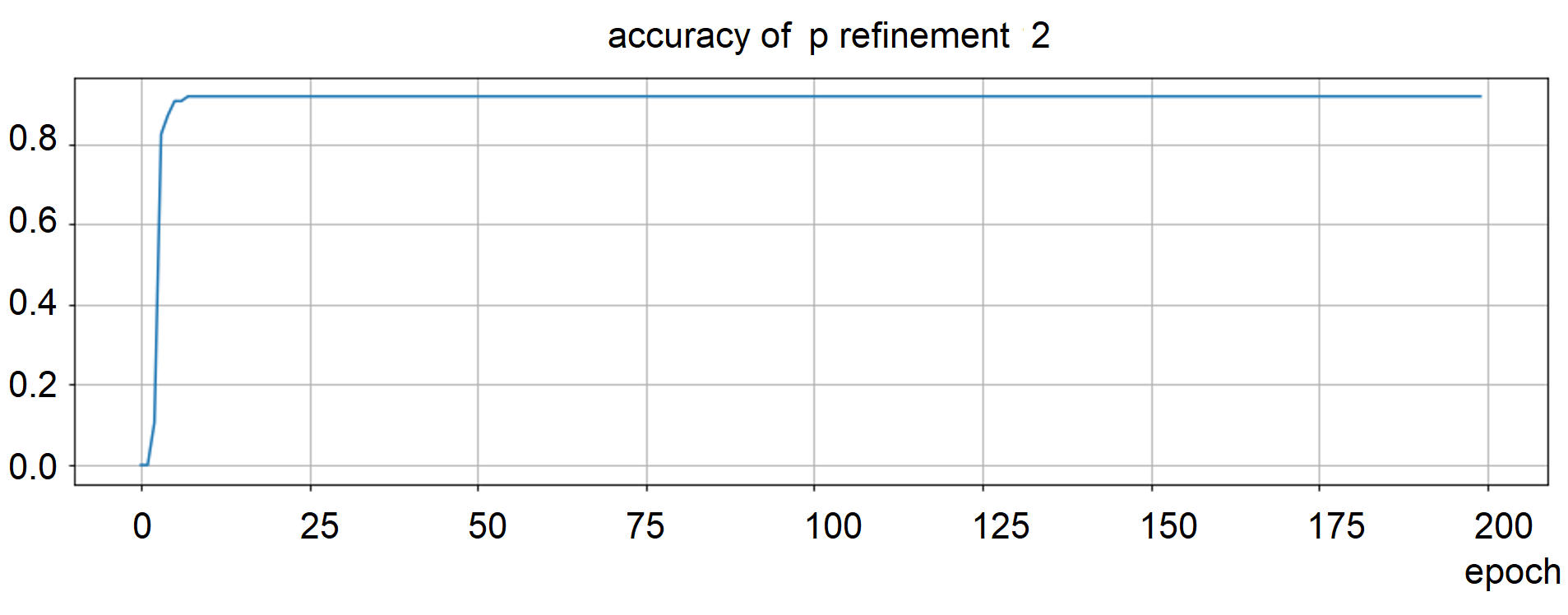}\\
    \includegraphics[width=0.49\textwidth]{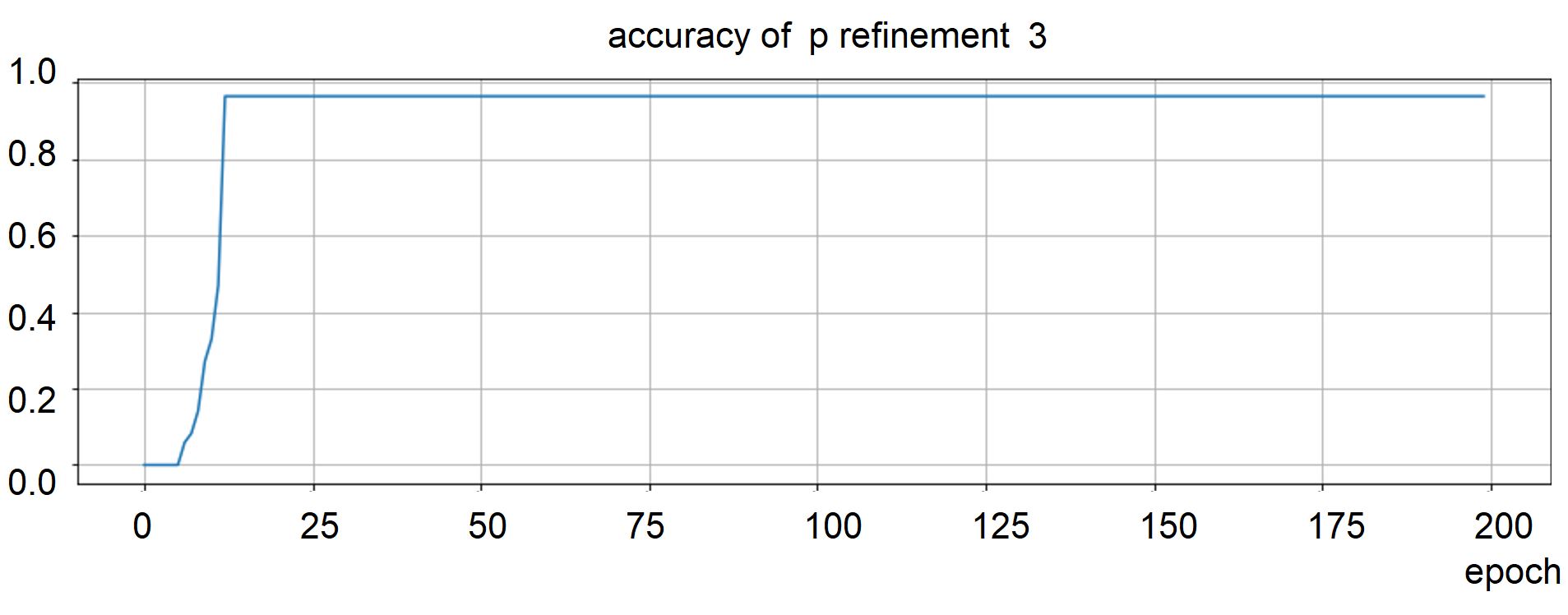}    \includegraphics[width=0.49\textwidth]{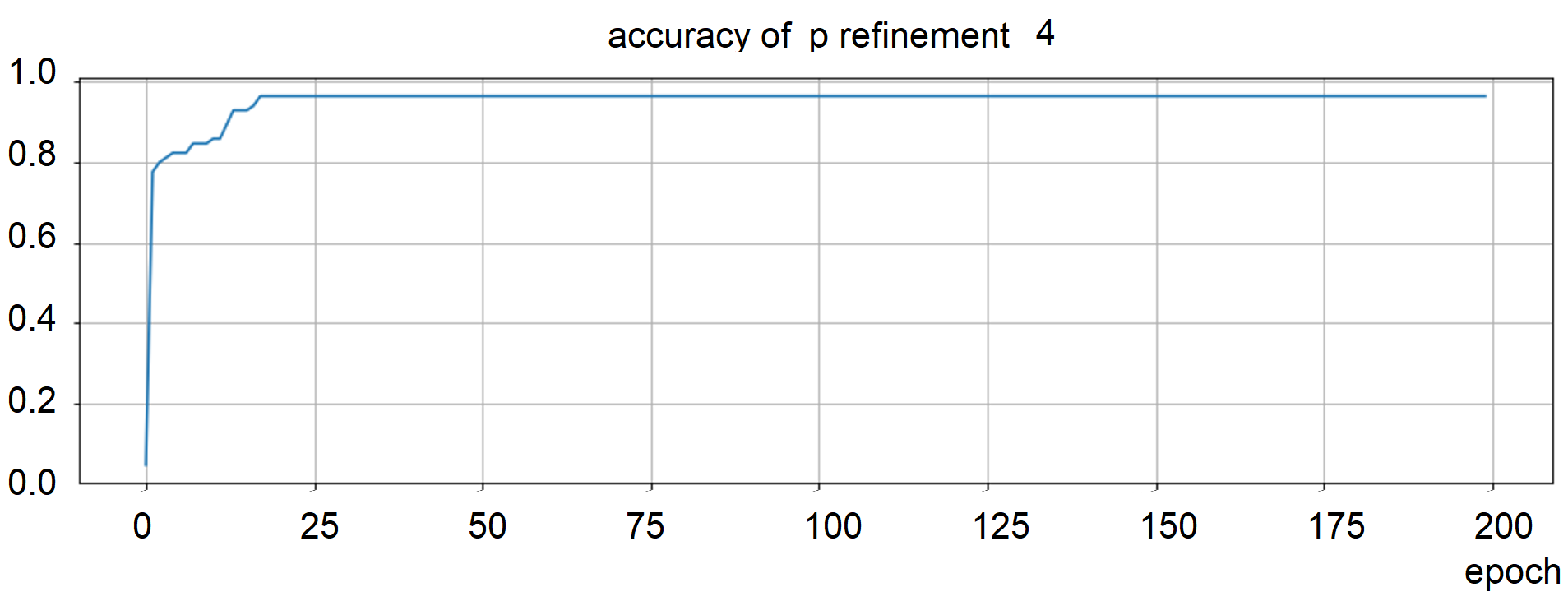}\\
    \includegraphics[width=0.49\textwidth]{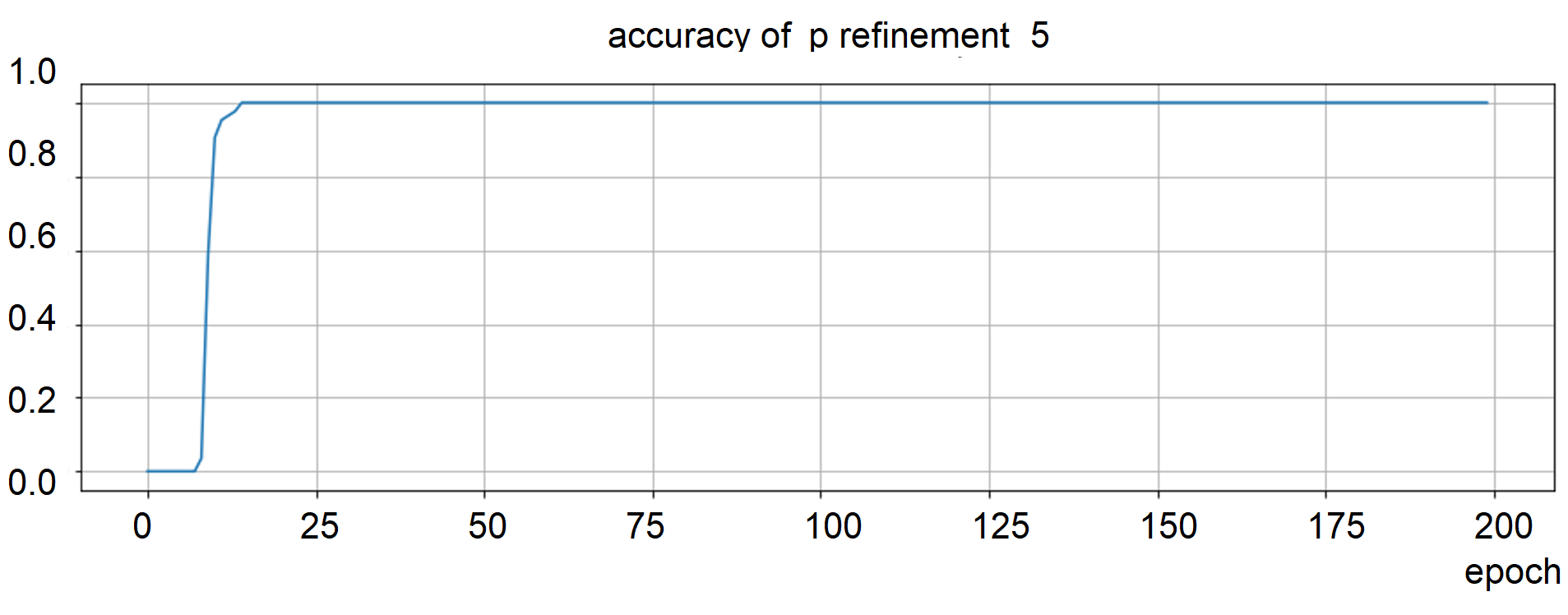}    \includegraphics[width=0.49\textwidth]{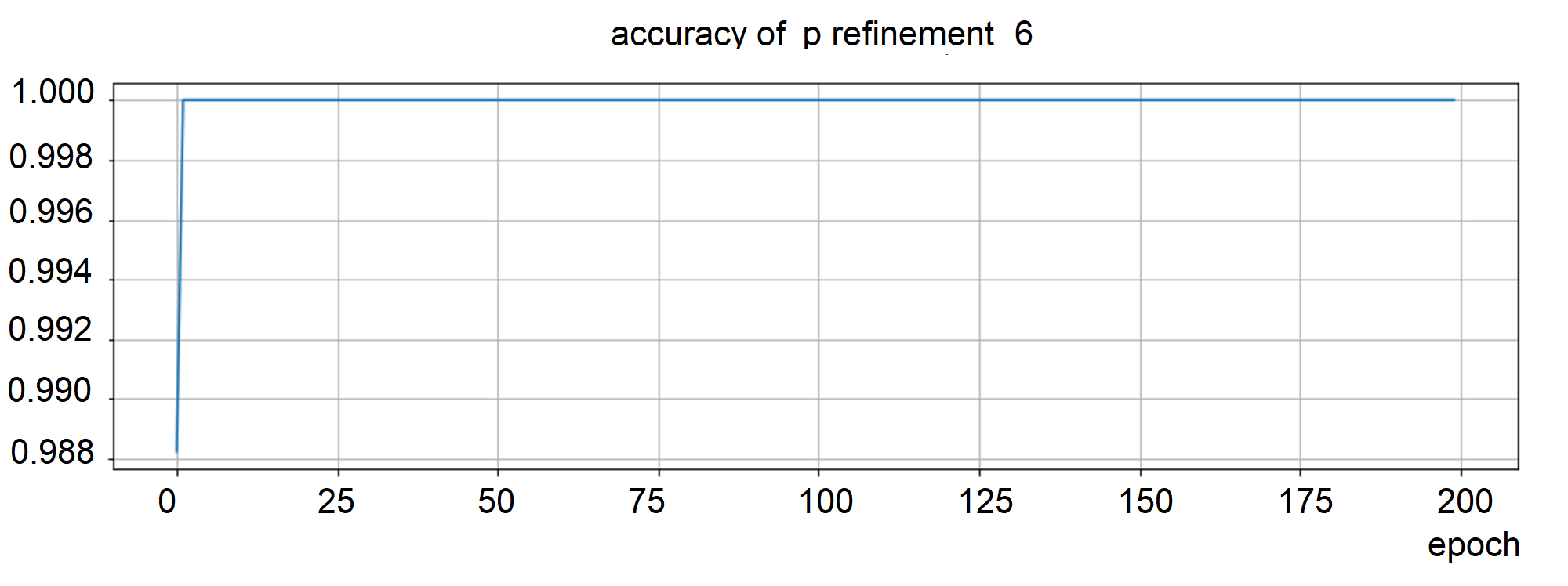}\\
    \includegraphics[width=0.49\textwidth]{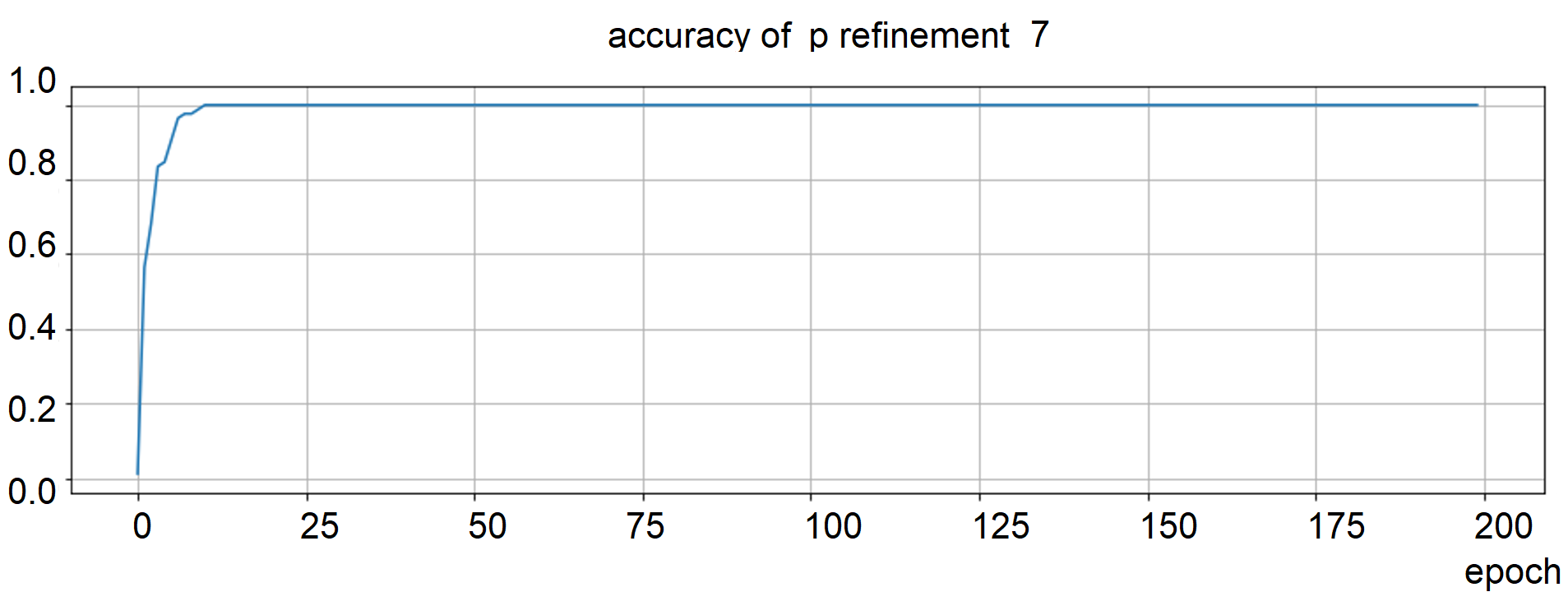}    \includegraphics[width=0.49\textwidth]{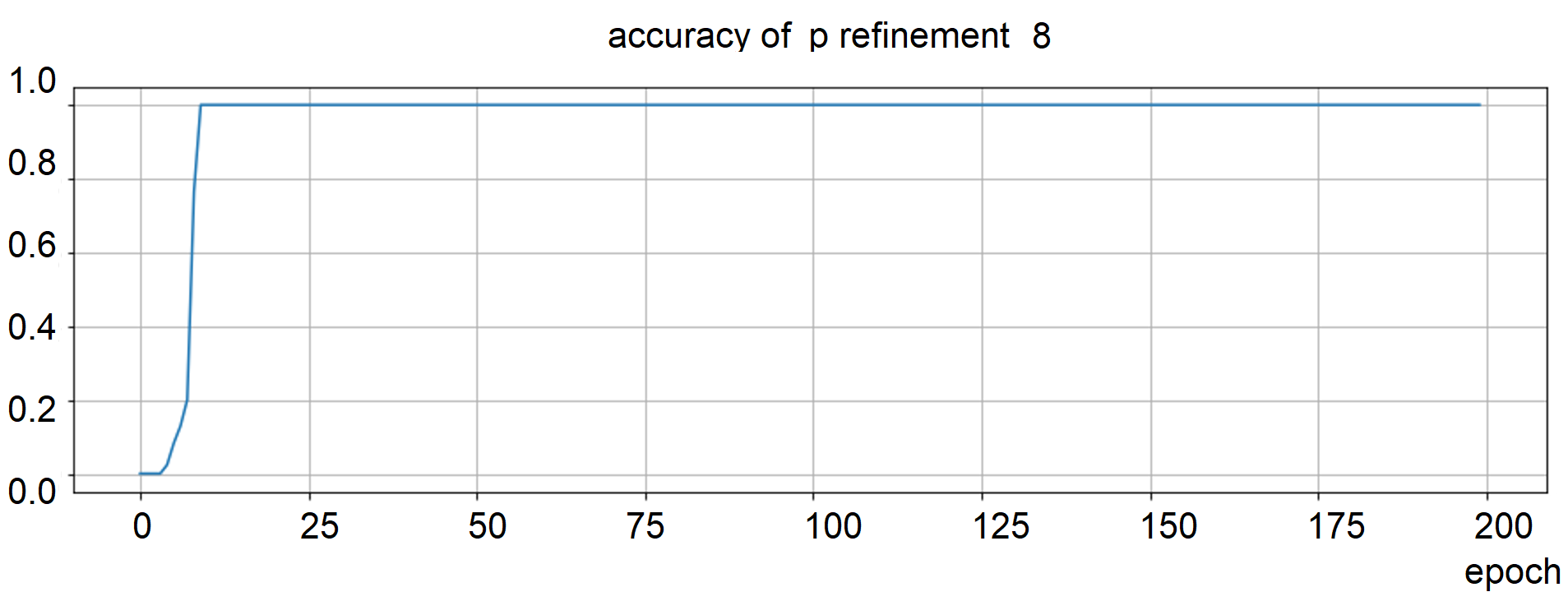}    
    \caption{The training procedure in 3D. The accuracy of the $h$ refinement. The accuracy of $p$ refinements for sons 1-8. }
    \label{fig:training3D}
\end{figure}

Figure \ref{fig:training3D} illustrates the training procedure in 3D. The dataset is partitioned into 80\% of training and 20\% validation sets. During each epoch, we train the DNN using the training dataset. After each training epoch, the DNN is tested against the validation subset. The plots in Figure \ref{fig:training3D} denote the percentage of the correctly predicted $h$ refinements and $p$ refinements taking into account all validation datasets.
Accuracy equal to 1.0 means the DNN makes 100 percent of incorrect decisions concerning the refinements.

The training in the 3D case is based on the 1132 samples generated during the execution of the twelve iterations of the self-adaptive $hp$-FEM algorithm on the Fichera model problem. The DNN has learned well the locations of the singularities, so it can propose optimal $h$ refinements along the three edges of the Fichera corner. The DNN has also learned how to distribute the polynomial orders of approximation. As we can read from the numerical results section, this enabled to continue with the exponential convergence.

\begin{figure}
    \centering
    \includegraphics[width=0.6\textwidth]{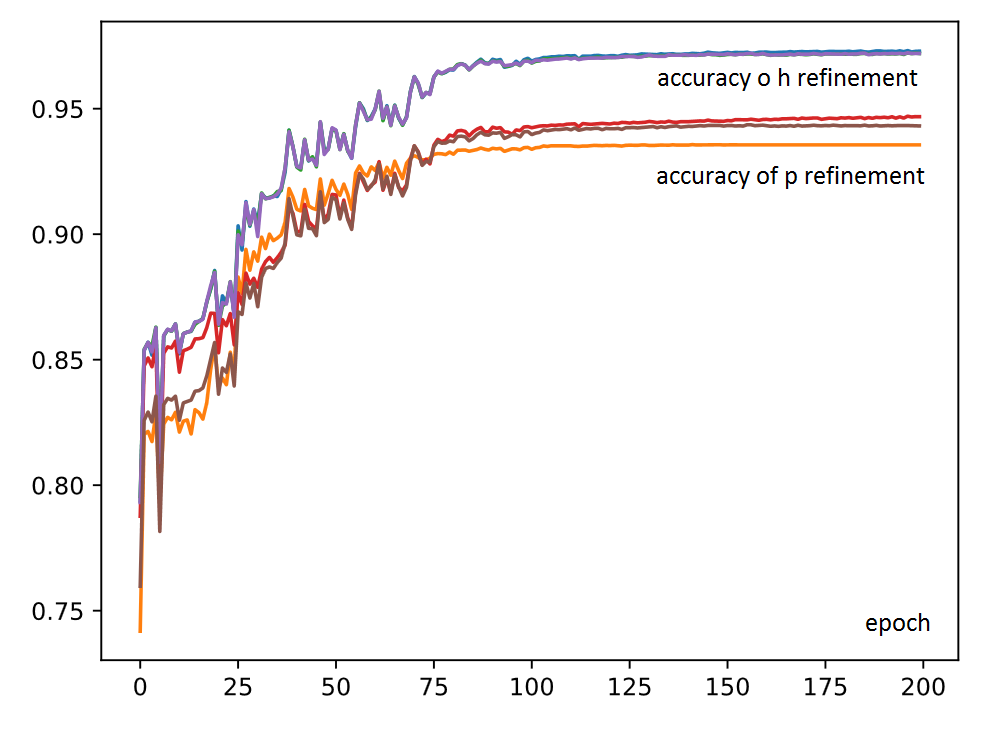}
    
    \caption{The training procedure in two-dimensions. }
    
\label{fig:training2D}
\end{figure}

Figure \ref{fig:training2D} illustrates the training procedure in two dimensions.
The dataset partitioning and the training procedure are similar to those in 3D. The plots in Figure \ref{fig:training2D} denote the percentage of the correctly predicted $h$ refinements and $p$ refinements considering all validation datasets.
Accuracy of $h$ refinement equal to 0.97 means the DNN makes 97 percent of correct decisions concerning the $h$ refinements. Accuracy equal to 0.95 means the DNN makes 5 percent of incorrect decisions concerning the $p$ refinements.

The training in the 2D case is based on the 10,000 samples generated during the 40 iterations of the self-adaptive $hp$-FEM algorithm on the L-shaped domain model problem.
In this case, learning how to $h$ refine means learning the location of the point singularity at point (0,0).
The DNN has learned the locations of this point singularity well so that it can propose the optimal $h$ refinements there. The DNN has also learned well the distribution of the polynomial orders of approximation from the 10,000 samples, with around 95 percent of correct decisions.
As discussed in the numerical results section, this enabled us to maintain exponential convergence.


\section{Conclusions}
We showed that an artificial expert can learn the optimal $hp$ refinement patterns for a given computational problem.
We proposed the Deep Neural Network driven $hp$-FEM algorithm (Algorithm 3), guiding the quasi-optimal refinements for the $hp$ finite element method.
The algorithm can be trained using the self-adaptive $hp$-FEM for a given problem. The self-adaptive $hp$-FEM uses an expensive direct solver to guide the optimal refinements. When we run out of resources, we can turn off the self-adaptive $hp$-FEM, collect the decisions about the optimal $hp$ refinements, and train the DNN. We can continue later with quasi-optimal refinements as guided by the DNN-driven $hp$-FEM. Instead of calling an expensive solver, we can ask the DNN expert to propose the $hp$ refinements. Our method has been verified 2D and 3D dimensions using the model L-shaped domain and the Fichera problems. In particular, we showed that we could obtain a high-quality DNN expert by using only element locations and $p$ refinement patterns.

\section*{Acknowledgments}
The European   Union's Horizon 2020 Research and Innovation Program of the Marie Sk\l{}odowska-Curie grant agreement No. 777778, MATHROCKs.
Research project partly supported by program "Excellence initiative – research university" for the University of Science and Technology.

\end{document}